\documentclass[11pt,a4paper]{article}
\usepackage[ansinew]{inputenc}
\usepackage[T1]{fontenc}

\usepackage{vmargin}
\usepackage{amsmath}
\usepackage{amsthm}
\usepackage{amssymb}
\usepackage{amsopn}
\usepackage{enumerate}
\usepackage{color}
\usepackage{esint}
\usepackage{mathtools}
\usepackage{hyperref}

\DeclareMathOperator{\Id}{Id}
\DeclareMathOperator{\dist}{dist}
\DeclareMathOperator{\diam}{diam}
\DeclareMathOperator{\co}{co}
\DeclareMathOperator{\spano}{span}
\DeclareMathOperator{\supp}{supp}

\theoremstyle{plain}
\newtheorem{thm}{Theorem}[section]
\newtheorem{prop}[thm]{Proposition}
\newtheorem{cor}[thm]{Corollary}
\newtheorem{lem}[thm]{Lemma} 

\theoremstyle{definition}
\newtheorem{defi}[thm]{Definition}

\newtheorem{ass}{Assumption}

\theoremstyle{remark}
\newtheorem{rem}[thm]{Remark}

\newcommand{\R}{\mathbb{R}}
\newcommand{\Q}{\mathbb{Q}}
\newcommand{\Z}{\mathbb{Z}}

\newcommand{\N}{\mathbb{N}}
\newcommand{\norm}[1]{\left\lVert#1\right\rVert}
\newcommand{\abs}[1]{\left\lvert#1\right\rvert}
\newcommand{\wto}{\rightharpoonup}

\newcommand{\beq}{\begin{equation}}
\newcommand{\eeq}{\end{equation}}
\newcommand{\bp}{\begin{proof}}
\newcommand{\ep}{\end{proof}}
\newcommand{\al}[1]{\begin{align*}#1\end{align*}}

\newcommand{\css}{\subset\subset}

\numberwithin{equation}{section}
\newcommand{\la}{\mathcal{L}}
\newcommand{\e}{\varepsilon}

\begin{document}

\begin{center}
\begin{Large}
On the passage from atomistic systems \\ to nonlinear elasticity theory
\end{Large}
\\[0.5cm]

Julian Braun\footnote{Institut f{\"u}r Mathematik, Universit{\"a}t Augsburg, Germany, {\tt julian.braun@math.uni-augsburg.de}} and Bernd Schmidt\footnote{Institut f{\"u}r Mathematik, Universit{\"a}t Augsburg, Germany, {\tt bernd.schmidt@math.uni-augsburg.de}}
\\[0.5cm]

\today
\\[1cm]
\end{center}

\begin{abstract}
We derive continuum limits of atomistic models in the realm of nonlinear elasticity theory rigorously as the interatomic distances tend to zero. In particular we obtain an integral functional acting on the deformation gradient in the continuum theory which depends on the underlying atomistic interaction potentials and the lattice geometry. The interaction potentials to which our theory applies are general finite range models on multilattices which in particular can also account for multi-pole interactions and bond-angle dependent contributions. Furthermore, we discuss the applicability of the Cauchy-Born rule. Our class of limiting energy densities consists of general quasiconvex functions and the class of linearized limiting energies consistent with the Cauchy-Born rule consists of general quadratic forms not restricted by the Cauchy relations.
\end{abstract}

\section{Introduction}

The main aim of this work is to provide a rigorous derivation of nonlinear elasticity functionals from atomistic models. The investigation of such discrete-to-continuum limits has been an active area of research in continuum mechanics over the last years in particular for, but not limited to, elastic interactions. For a recent account on this line of research and a summary of the related literature we refer to the survey article \cite{BLL:07} by Blanc, LeBris and Lions.

Classically, the stored energy density in elasticity theory is derived from atomistic models by applying the Cauchy-Born rule: Given a macroscopic deformation $y$ of the elastic body, one assumes that microscopically near every material point $x$, all the atoms deform by just following the macroscopic deformation gradient $F = \nabla y(x)$. Inserting this ansatz into the atomistic potentials then leads to a continuum stored energy density $W$ as a function of $F \in \R^{3 \times 3}$. Assuming validity of the Cauchy-Born rule, very general and even quantum mechanical interactions have, e.g., been investigated by Blanc, LeBris and Lions in \cite{BLL:02}. A priori, however, it is not clear if the Cauchy-Born hypothesis does hold true. 

For a two-dimensional mass-spring model, it has been shown by Friesecke and Theil in \cite{friesecketheil} that the Cauchy-Born rule does indeed hold true for small strains, while it in general fails for large strains. This result has then be generalized to a wider class of discrete models and more than two dimensions by Conti, Dolzmann, Kirchheim and M{\"u}ller in \cite{CDKM}. 

A fundamental contribution towards a rigorous derivation of continuum limits in elasticity has been made by Alicandro and Cicalese in \cite{alicandrocicalese}, where they prove a general integral representation result for continuum limits of atomistic pair interaction potentials. It is our main aim, departing from this result, to derive a continuum theory for more general interaction potentials which, in particular, can also incorporate bond-angle dependent potentials. Such an extension is desirable in applications, as many atomistic models such as, e.g., the Stillinger-Weber potential, cannot be written as a pure pair potential. In fact, the class of potentials our theory applies to is rich enough to model any continuum energy density, even if the Cauchy-relations are violated. See below for more details. For small strains, such general models have been derived by the second author in \cite{schmidtlinelast}. A first step in this direction in the nonlinear regime has recently been provided through the analysis of a special class of nearest neighbor three point interactions on a two-dimensional square lattice by Meunier, Pantz and Raoult, see \cite{MPR}, and for special interactions subordinate to specific triangulations by Alicandro, Cicalese and Gloria \cite{alicandrocicalesegloria}.
 
Our limiting density will be described in terms of a sequence of cell problems. This is related to the homogenization results of Braides \cite{braides-hom} and M{\"u}ller \cite{mueller-hom} of nonconvex integral functionals. In the setting of a discrete-to-continuum limit for thin films in the membrane energy regime, such a limiting cell formula has also been obtained by the second author in \cite{schmidt-membrane}. Indeed, following the localization method (cf., e.g, \cite{dalmasogamma}), also from a technical point of view our reasoning is related to \cite{braides-hom}.
 
The atomistic systems we will consider are `generalized mass spring models' as in \cite{CDKM}. For such models, a rigorous simultaneous discrete-to-continuum and nonlinear-to-linear limit has been obtained by the second author in \cite{schmidtlinelast}. The present paper now extends these results in the purely nonlinear setting. While in the linearized regime, as shown in \cite{schmidtlinelast}, the limiting energy is indeed given by the Cauchy-Born energy, in the nonlinear setting this cannot be true in general, cf.\ \cite{friesecketheil}. Nevertheless, we will obtain a definite stored energy density in the limit, which under appropriate conditions will be equal to the Cauchy-Born density for small (but finite) strains. 
 
In order to prove our main representation result we resort to abstract compactness properties of $\Gamma$-limits and integral representation results for functionals on Sobolev spaces and thus follow the scheme set forth in \cite{alicandrocicalese}, which is dictated by verifying the hypotheses of that abstract approach by the localization method. A few of the arguments in this proof can be used with only minor adjustments. There are, however, some major differences as compared to the pair interaction case treated by Alicandro and Cicalese.
While these authors use slicing arguments in order to obtain energy estimates on the usual $d \times d$ deformation gradients in the direction of interacting pairs, we will have to estimate the much higher dimensional $d \times 2^d$ discrete deformation gradients. In fact, as in general our discrete energies cannot be recovered by slicing techniques, we will instead work with very carefully chosen interpolations of the discrete deformations which encode the full discrete gradient on lattice cells. 

More specifically, if ${\cal L} = A \Z^d$ is some Bravais lattice, $\Omega \subset \R^d$ a bounded open set with Lipschitz boundary that will be viewed as the `macroscopic' domain occupied by the elastic body, whose atoms are at positions $\e {\cal L} \cap \Omega$, we assume that the energy of a deformation $y : \e {\cal L} \cap \Omega \to \R^d$ can be written in the form 
$$ F_{\e}(y, \Omega) = \e^d \sum\limits_{x\in (\la_\e'(\Omega))^\circ} W_{\rm cell}(\bar{\nabla}y(x)) , $$
where $x$ runs over all midpoints of elementary lattice cells of $\e {\cal L}$ inside $\Omega$. Here $\bar{\nabla}y(x)$ is the discrete gradient of $y$ on the corresponding cell $Q$ which encodes all relative displacements of atoms lying on the corners of $Q$. (See Section \ref{sec:prelim} for precise definitions.) $\e$ is the small parameter in the system measuring the typical interatomic distance and tending to zero eventually in the continuum limit. The rescaling by $\e^d$ is introduced in order to pass from units of finite energy per atom to units of finite energy per unit volume. 

In fact, our analysis is not restricted to interactions within unit lattice cells, but also applies to general finite range interactions. In such models, the energy is still given as the sum over unit lattice cells $(\la_\e'(\Omega))^\circ$, but the cell energy now depends on the discrete deformation gradient $\bar{\nabla}y(x) \in \R^{d \times N}$ on a larger `super-cell': 
\[ F_{\e}(y, \Omega) = \e^d \sum\limits_{x\in (\la_\e'(\Omega))^\circ} W_{\rm super-cell}(\bar{\nabla}y(x)), \]
where the definition of the lattice interior $(\la_\e'(\Omega))^\circ$ is suitably adjusted so that only lattice points within $\Omega$ may interact.

Another complication arises when extending our results to general finite-range interactions on multi-lattices. For such systems, the discrete gradients are augmented with additional internal variables describing the relative shifts of the underlying single lattice. With the help of a mixed Sobolev/Lebesgue space representation theorem we are then led to a boundary value/mean value cell formula for the limiting energy density. This cell formula is in fact related to the cell formula derived in \cite{schmidt-membrane} for thin membranes where internal variables measure relative shifts of the thin film's layers.
On multi-lattices $\e (\{0, s_1, \ldots, s_m\} + {\cal L})$ with $m$ shift vectors $s_1, \ldots, s_m \in \R^d$ the general discrete energy functional then reads 
$$ F_{\e}(y, s, \Omega) = \e^d \sum\limits_{x\in (\la_\e'(\Omega))^\circ} W_{\rm super-cell}(\bar{\nabla}y(x), s(x)), $$
where $y : \e {\cal L} \cap \Omega \to \R^d$, $s : \e {\cal L} \cap \Omega \to \R^{d\times m}$.
For notational convenience we will restrict to simple unit cell interactions on a Bravais lattice for the largest part of the paper and only comment on the necessary modifications in the more general case at the end of Section \ref{sec:hom}. 

Our main results are summarized in the following theorems. The necessary assumptions Assumptions \ref{ass:coercive}, \ref{ass:growth} and \ref{ass:CB} on the cell energy are specified in Section \ref{sec:prelim}. (Assumptions \ref{ass:coercive} and \ref{ass:growth} are nothing but standard $p$-growth assumptions on $W_{\rm cell}$.)
\begin{thm}[$\Gamma$-convergence]\label{thm:homogenization}
	Suppose Assumptions \ref{ass:coercive} and \ref{ass:growth} are true. Then $F_\e(\cdot,\Omega)$ $\Gamma(L^p(\Omega;\R^d))$- and $\Gamma(L^p_{loc}(\Omega;\R^d)/\R)$-converges to the functional $F$, defined by
	\[
	F(y) =
	\begin{dcases}
		\int\limits_\Omega W_{\rm cont}(\nabla y (x)) \,dx, &\text{ if } y\in W^{1,p}(\Omega;\R^d),\\
		\infty &\text{ otherwise,}
	\end{dcases}
	\]
	where the continuum density $W_{\rm cont}\colon \R^{d\times d}\to [0,\infty)$ is given in terms of $W_{\rm cell}$ by
	\[
		W_{\rm cont}(M) = \frac{1}{\abs{\det A}} \lim\limits_{N\to\infty} \frac{1}{N^d} \inf\left\{ \sum\limits_{x\in (\la_1'(A(0,N)^d))^\circ} W_{\rm cell}(\bar{\nabla}y(x)) \colon y \in \mathcal{B}_1(A(0,N)^d,y_M) \right\}.
	\]
\end{thm}
Here $\mathcal{B}_1(A(0,N)^d,y_M)$ is the space of lattice deformations of $\la_1 \cap A(0,N)^d$ with linear boundary conditions $M$ on $\partial \la_1(A(0, 1)^d)$, cf.\ Section \ref{sec:bvp}. 

As in non-convex homogenization (see \cite{mueller-hom}), in the representation result for $W_{\rm cont}$ it is nesessary to minimize $W_{\rm cell}$ over larger and larger cubes and the limit is in general not obtained for finite $N$. A simple 2d example for this effect is given by a square lattice where nearest neighbor atoms interact via a harmonic spring potential: $F_{\e}(y) = \frac{1}{2} \sum_{|x - x'|=\e}(|y(x) - y(x')| - \e)^2$ (which can be written in the above form). This is a simplified version of the example discussed in \cite[Sect.~4.4]{schmidt-membrane}. The arguments sketched there, which amount to considering deformations 
\[ y(x_1, x_2) = M(x_1, x_2) + \sigma_1(x_1) + \sigma_2(x_2) \] 
for compressive boundary conditions, with suitable $2$-periodic functions $\sigma_1$ and $\sigma_2$ like
\[\sigma_i(z) = \frac{1}{2}(-1)^z\sqrt{\frac{1}{m_{1i}^2+m_{2i}^2}-1}
\begin{pmatrix}
-m_{2i}\\
m_{1i}
\end{pmatrix},  \]
$i = 1,2$, and suitably modified on the boundary, show that 
\begin{align*}
  \frac{1}{N^2} \inf\left\{ \sum\limits_{x\in (\la_1'((0,N)^d))^\circ} 
  W_{\rm cell}(\bar{\nabla}y(x)) \colon y \in \mathcal{B}_1((0,N)^d,y_M) \right\}
\end{align*}
converges to $(\max\{0, |m_{\cdot 1}| - 1\})^2 + (\max\{0, |m_{\cdot 2}| - 1\})^2$ with error bound $O(N^{-1})$, where $m_{\cdot j}$ denotes the $j$th column of $M$. Evaluating, however, the energy of the bonds on and close to the boundary it is easily seen that, for compressive boundary conditions $|m_{\cdot 1}| < 1$ or $|m_{\cdot 2}| < 1$, this limiting energy is always over-estimated by a constant times $N^{-1}$.

\begin{thm}[Compactness]\label{thm:bd-En-seq}
Under the assumptions of Theorem \ref{thm:homogenization}, if $y_{\e}$ is a sequence with equibounded energies $F_{\e}(y_{\e}, \Omega)$ and $\Omega$ is connected, then there exist a sequence $\e_k \to 0$ and $y \in W^{1, p}(\Omega; \R^d)$ such that $y_{\e_k} \to y$ in $L^p_{loc}(\Omega;\R^d)/\R$. 
\end{thm}
Of course, if $\Omega$ is not connected, one has compactness in $L^p_{loc}$ up to translation on every connected component.
	
Analogous results hold true under boundary conditions $g\in W^{1,\infty}(\R^d;\R^d)$. Let $F^g$ and $F_\e^g$ denote the functionals obtained from $F$ and $F_{\e}$, respectively, with values set to infinity if the boundary conditions are not met, cf.\ Section \ref{sec:bvp}. 
\begin{thm}[$\Gamma$-convergence]\label{thm:homogenization-bdry}
	If Assumptions \ref{ass:coercive} and \ref{ass:growth} are true, then $F_\e^g(\cdot,\Omega)$ $\Gamma(L^p(\Omega;\R^d))$-converges to the functional $F^g$.
\end{thm}
\begin{thm}[Compactness]\label{thm:bd-En-seq-bdry}
Under the assumptions of Theorem \ref{thm:homogenization-bdry}, if $y_{\e}$ is a sequence with equibounded energies $F_{\e}^g(y_{\e}, \Omega)$, then there exist a sequence $\e_k \to 0$ and $y \in W^{1, p}(\Omega; \R^d)$ with $y = g$ on $\partial \Omega$ such that $y_{\e_k} \to y$ in $L^p(\Omega;\R^d)$. 
\end{thm}

A standard argument then yields that almost minimizers of $F^g_{\e}(\cdot, \Omega)$ converge to minimizers of $F^g(\cdot, \Omega)$ and almost minimizers of $F_{\e}(\cdot, \Omega)$ up to translation converge to minimizers of $F(\cdot, \Omega)$, more precisely:
\begin{cor}[Convergence of almost minimizers]\label{cor:conv-almost-min}
	Suppose Assumptions \ref{ass:coercive} and \ref{ass:growth} are true. Then every sequence of almost minimizers of $F_{\e}(\cdot, \Omega)$ for connected $\Omega$ is compact in $L^p_{loc}(\Omega;\R^d)/\R$ and every limit is a minimizer of $F$, while every sequence of almost minimizers of $F^g_{\e}(\cdot, \Omega)$ is compact in $L^p(\Omega;\R^d)$ and every limit is a minimizer of $F^g$.
\end{cor}

It is not hard to include body forces in the energy expression as these will only be continuous perturbations of the energy functional which converge uniformly on bounded sets and thus preserve $\Gamma$-convergence by general theory.

We also remark that the point why the theory can be adapted to the case of general finite range interactions, is that in this case $W_{\rm cell}$, while naturally still bounded from above by the discrete gradient through Assumption \ref{ass:growth}, from below has to be bounded only in terms of the discrete gradient on the unit cell. See Section \ref{sec:hom} for details.
For general finite range interactions on multi-lattices we state here only the analogue of Theroem \ref{thm:homogenization}, as in the cell formula there are now additional internal variables that need to be taken into account. Theorems \ref{thm:bd-En-seq}, \ref{thm:homogenization-bdry} and \ref{thm:bd-En-seq-bdry} and Corollary \ref{cor:conv-almost-min} extend in a straightforward manner. 
\begin{thm}[$\Gamma$-convergence]\label{thm:homogenization-multicell}
	Suppose $W_{\rm super-cell}$ satisfies the growth assumptions (as stated in Section \ref{sec:hom}). Then $F_\e(\cdot,\cdot,\Omega)$ $\Gamma(L^p(\Omega;\R^d) \times \mbox{w-}L^q(\Omega;\R^{d \times m}))$-converges to the functional $F$, defined by
\begin{align*}
	F(y,s) =
	\begin{dcases}
		\int\limits_\Omega W_{\rm cont}(\nabla y (x), s(x)) \,dx, 
		&\text{ if } y\in W^{1,p}(\Omega;\R^d), \\
		\infty &\text{ otherwise,}
	\end{dcases}
\end{align*}
	where the continuum density $W_{\rm cont}\colon \R^{d\times d} \times \R^{d \times m}\to [0,\infty)$ is given in terms of $W_{\rm super-cell}$ by
\begin{align*}
  W_{\rm cont}(M, s_0) 
  = \frac{1}{\abs{\det A}} \lim\limits_{N\to\infty} \frac{1}{N^d} \inf 
     \Bigg\{ \sum\limits_{x\in (\la_1'(A(0,N)^d))^\circ} W_{\rm super-cell}(\bar{\nabla}y(x), s(x)) \colon &\\ 
     (y,s) \in \mathcal{B}_1(A(0,N)^d,y_M,s_0) \Bigg\}&.
\end{align*}
\end{thm}
Here $\mathcal{B}_1(A(0,N)^d,y_M,s_0)$ is the space of lattice deformations of $\la_1 \cap A(0,N)^d$ with linear boundary conditions $M$ on $\partial \la_1(A(0, 1)^d)$ for $y$ and average $s_0$ for $s$, cf.\ Section \ref{sec:hom}.

As macroscopic deformations are solely given in terms of a deformation mapping $y\in W^{1,p}(\Omega;\R^d)$, we are also interested in the effective macroscopic energy density obtained by minimizing out the internal variables $s$:
\begin{thm} \label{thm:mins}
	For every $y\in W^{1,p}(\Omega;\R^d)$ we have
	\[ \min\limits_{s\in L^q(\Omega;\R^{d\times m})} \int\limits_\Omega W_{\rm cont}(\nabla y(x),s(x))\, dx = \int\limits_\Omega \min\limits_{s\in \R^{d\times m}} W_{\rm cont}(\nabla y (x),s)\, dx.\]
Moreover, $F^{s-\min}_\e(\cdot,\Omega) := \inf\limits_{s \in L^q} F_{\e}(\cdot,s,\Omega)$ $\Gamma(L^p(\Omega;\R^d))$-converges to the functional $F^{s-\min}$, defined by
\begin{align*}
	F^{s-\min}(y) =
	\begin{dcases}
		\int\limits_\Omega \min\limits_{s\in \R^{d\times m}} W_{\rm cont}(\nabla y (x), s) \,dx, 
		&\text{ if } y\in W^{1,p}(\Omega;\R^d), \\
		\infty &\text{ otherwise,}
	\end{dcases}
\end{align*}
\end{thm}

Returning to our basic setting on a Bravais lattice, under an additional assumption, we can calculate the limiting density for small strains explicitly by the Cauchy-Born rule: 
\begin{thm}\label{thm:CB}
	In addition to Assumptions \ref{ass:coercive} and \ref{ass:growth} suppose that $W_{\rm cell}$ satisfies Assumption \ref{ass:CB}. Then there is a neighborhood ${\cal U}$ of $SO(d)$ such that $W_{\rm cont}$ is given by 
	\[  
	   W_{\rm cont}(M) = W_{\rm CB}(M) := \frac{1}{\abs{\det A}} W_{\rm cell}(MZ). 
	\]
	for all $M \in {\cal U}$.

\end{thm}
Here $Z \in \R^{d \times 2^d}$ is a `discrete identity matrix', see Section \ref{sec:prelim} for details. 

As $W_{\rm cont}$ arises as the energy density of a $\Gamma$-limit it has to be quasiconvex (cf.\ Section \ref{sec:prelim} for the definition of these concepts). The next proposition shows that our class of atomistic interactions is rich enough to model any quasiconvex energy density in the continuum limit.

\begin{prop}\label{prop:all-qc}
Suppose $V : \R^{d \times d} \to \R$ is quasiconvex with standard $p$-growth 
\[ c \abs{M}^p - c' \le V(M) \le c'' (\abs{M}^p +1) \]
for some constants $c, c', c'' > 0$ and all $M \in \R^{d \times d}$. Then there exists a cell energy $W_{\rm cell}$ satisfying Assumptions \ref{ass:coercive} and \ref{ass:growth} such that $W_{\rm cont} = V$. 
\end{prop}

We remark that, by way of contrast, a restriction to pair interaction models will only lead to a restricted class of limiting continuum energies, as can be quantified in terms of the so-called Cauchy relations: If the Cauchy-Born rule applies (e.g., due to Assumption \ref{ass:CB}), an atomistic interaction energy
\[ E(y)=\e^d\sum\limits_{x,x'\in\la_{\e} \cap \Omega \atop x\neq x'} V_{\frac{|x-x'|}{\e}}\left(\frac{\abs{y(x)-y(x')}}{\e}\right) \] 
yields the continuum density
\[ W_{\rm CB}(M)=\frac{1}{\abs{\det A}}\sum\limits_{x\in\la \atop x \neq 0} V_{|x|}(\abs{Mx}). \]
Assuming $V_{|x|}$ is smooth and, for large $|x|$, sufficiently rapidly decreasing a direct calculation yields
\[ D^2 W_{\rm CB}(\Id)(M, M) 
=  \sum\limits_{i,j,k,l=1}^d c_{ijkl} m_{ij} m_{kl}, \]
where the elastic constants $c_{ijkl}$ are given by 
\[ c_{ijkl} 
= \frac{1}{\abs{\det A}}\sum\limits_{x\in\la \atop x \neq 0} V_{|x|}''(\abs{x}) \frac{x_i x_j x_k x_l}{\abs{x}^2} 
+ V_{|x|}'(\abs{x}) \bigg( \frac{x_j x_l \delta_{ki}}{\abs{x}} - \frac{x_i x_j x_k x_l}{\abs{x}^3} \bigg). \]
While the symmetry relations $c_{ijkl} = c_{klij}$ and $c_{ijkl} = c_{jikl}$ naturally follow from the symmetry of the Hessian $D^2 W_{\rm CB}(\Id)$ and frame indifference of $W_{\rm CB}$, the particular form of $W_{\rm CB}$ in addition gives $c_{ijkl} = c_{ilkj} = c_{kjil}$ for every $i,j,k,l$.

In the 3-dimensional setting of elasticity theory these additional relations lower the dimension of admissible elasticity tensors from $21$ to $15$ (symmetric in all indices) and so can be written as $6$ equations, the Cauchy-relations, namely
\al{ c_{1122} = c_{1212},\quad
c_{2233} = c_{2323},\quad
c_{3311} = c_{3131}, \\
c_{1123} = c_{1213},\quad
c_{2231} = c_{2321},\quad
c_{3312} = c_{3132}. }

The question whether in elasticity theory the Cauchy-relations hold true (rari-constant theory) or fail (multi-constant theory) had been under discussion for quite some time in physics and was finally decided by experimental data in favour of the multi-constant theory (for some experimental data and further physical considerations cf.\ \cite{haussuehl}). This means, in particular, that the interaction in a lattice is a complex multibody interaction which cannot be reduced to pair-potentials.
Our model in this paper using general cell energies is not limited by the Cauchy-relations: 

\begin{prop}\label{prop:all-quad}
Suppose $Q : \R^{d \times d} \to \R$ is a positive semidefinite quadratic form which is positive definite on the symmetric $d \times d$ matrices and vanishes on antisymmetric matrices. Then there exists $W_{\rm cell}$ satisfying Assumptions \ref{ass:coercive}, \ref{ass:growth} and \ref{ass:CB} such that 
\[ \frac{1}{2} D^2 W_{\rm CB}(\Id)(M,M) = Q(M). \]
\end{prop}

The paper is organized as follows. In Section \ref{sec:prelim} we first introduce the discrete model and review some basic facts on $\Gamma$-convergence and integral representations of functionals on Sobolev spaces. In Section \ref{sec:repr} we then proceed to state precisely and prove a general $\Gamma$-compactness and representation theorem. This in particular requires a number of technical preliminaries in order to investigate discrete deformations. A version of this representation result for boundary value problems is then provided in Section \ref{sec:bvp}. Finally, the limiting stored energy function is identified in Section \ref{sec:hom} through minimizing a sequence of cell problems, leading to the main discrete-to-continuum convergence result and the proofs of the results stated above.

\section{ The model and general preliminaries}\label{sec:prelim}

In this section we introduce the atomistic model and recall some general facts on $\Gamma$-convergence and integral representation results required by the localization method.

\subsection{The atomistic model}

Let $\la \subset \R^d$ be a Bravais lattice, i.e., there are linearly independent vectors $v_1, \dots, v_d $ such that \[\la = \{n_1 v_1 + \dots + n_d v_d \mid n_1, \dots, n_d \in \Z\}= A \Z^d,\]
if we set $A=(v_1, \dots, v_d)$.
The scaled lattices $\la_\e = \e \la$ partition $\R^d$ into the $\e$-cells $z + A[0,\e)^d$ ($z\in\la_\e$). Let $Q_\e(x)$ denote the $\e$-cell containing $x$. The centers of the cells are $\la_\e' = \la_\e + A(\frac{1}{2}, \dots, \frac{1}{2})$ and we denote by $\bar{x}$ the center of the cell containing $x$. These centers give a convenient labeling of the cells. 
Furthermore let $z_1, \dots, z_{2^d}$ be the points in $A\left\{-\frac{1}{2},\frac{1}{2}\right\}^d$ and $Z:=(z_1, \dots,z_{2^d}) \in \R^{d \times 2^d}$.

For a set $U\subset\R^d$ we define the following lattice subsets in the spirit of its closed hull, interior or boundary with respect to $\e\la'$ or its corners $\e\la$ by 
\[ \la_\e'(U) = \{ x \in \la_\e' \mid \overline{Q_\e(x)} \cap U \neq \emptyset \}, \quad \la_\e(U) = \la_\e'(U) + \e\{z_1, \dots, z_{2^d}\} ,\]
\[ (\la_\e'(U))^\circ = \{ x \in \la_\e' \mid \overline{Q_\e(x)} \subset U \}, \quad (\la_\e(U))^\circ = (\la_\e'(U))^\circ + \e\{z_1, \dots, z_{2^d}\} ,\]
\[\partial\la_\e'(U) = \la_\e'(U) \backslash (\la_\e'(U))^\circ, \quad \partial\la_\e(U) = \partial\la_\e'(U) + \e\{z_1, \dots, z_{2^d}\}.  \]
Furthermore let \[U^\e = \bigcup\limits_{\bar{x} \in \la_\e'(U)} Q_\e(\bar{x}),\quad U_\e = \bigcup\limits_{\bar{x} \in (\la_\e'(U))^\circ} Q_\e(\bar{x}).\]

A lattice deformation should be thought of as a mapping $\la_\e \cap \Omega \to \R^d$. Choosing a suitable extension (e.g., by $0$) and piecewise constant interpolation, we can and will assume that the lattice deformations $\mathcal{B}_{\e}(\Omega)$ are the functions $\Omega \to \R^d$, which are constant on every cell $Q_\e(x)$, $x\in \la_\e'(\Omega)$. (This will not change the energy, see below.) 

If we have a deformation $y \in \mathcal{B}_{\e}(\Omega)$ and $x\in \Omega_\e$, we set $y_i(x) = y(\bar{x} + \e z_i)$,
\[\bar{y}(x)=\frac{1}{2^d} \sum\limits_{i=1}^{2^d} y_i(x) \quad \text{and} \quad \bar{\nabla}y(x) = \frac{1}{\e} (y_1(x)-\bar{y}(x), \dots, y_{2^d}(x)-\bar{y}(x)). \]

Let $\mathcal{A}(U)$ be the set of all bounded open subsets of $U\subset\R^d$ and $\mathcal{A}_L(U)$ the set of all those, that have a Lipschitz boundary.
In the following, we will consider a set $\Omega\in\mathcal{A}_L(\R^d)$ and the energies $F_\e \colon L^p(\Omega;\R^d) \times \mathcal{A}(\Omega) \to [0,\infty]$ for some fixed $1<p<\infty$, defined by
\beq
 F_\e(y,U) = 
 \begin{dcases}
  \e^d \sum\limits_{x\in (\la_\e'(U))^\circ} W_{\rm cell}(\bar{\nabla}y(x)) & \text{if } y\in\mathcal{B}_\e(U), \\
  \infty  & \text{otherwise.}
 \end{dcases}
\eeq
In this definition the energy only depends on the values of $y$ in $(\la_\e(U))^\circ\subset \la_\e \cap U$. Of course, there can be some points in $\la_\e \cap U$ which we do not use at all, but this is negligible if we impose Dirichlet boundary conditions as we will do later on.

We make some assumptions on the cell energy $W_{\rm cell} \colon \R^{d \times 2^d} \to [0,\infty)$. Note, that a discrete gradient can take values precisely in the space
\[ V_0 = \left\{ F\in\R^{d\times 2^d} \colon \sum\limits_{j=1}^{2^d} a_{ij}=0,\text{ for every } i=1, \dotsc, d  \right\}.\]
Therefore, we are only interested in the values of $W_{\rm cell}$ on $V_0$.

\begin{ass} \label{ass:coercive}
There are $c,c'>0$ such that for every $F\in V_0$
 \[W_{\rm cell}(F) \geq c\abs{F}^p-c'.\]
\end{ass}
\begin{ass} \label{ass:growth}
There is a $c>0$ such that for every $F\in V_0$
 \[W_{\rm cell}(F) \leq c(\abs{F}^p+1).\]
\end{ass}

While these conditions are supposed to hold true for all our results, we also state a third assumption, which, if satisfied, allows for an application of the Cauchy-Born rule locally near $SO(d)$. The so-called Cauchy-Born energy density is defined by letting each atom follow the macroscopic gradient: 
$$ W_{\rm CB}(M) := \frac{1}{\abs{\det A}} W_{\rm cell}(MZ) $$ 
for $M \in \R^{d \times d}$. 
 
\begin{ass} \label{ass:CB}
\begin{enumerate}[(i)] 
\item $W_{\rm cell} : \R^{d \times 2^d} \to \R$ is invariant under translations and rotations, i.e.\ for $F \in \R^{d \times 2^d}$, 
\begin{align*} 
  W_{\rm cell} (R F + (c, \ldots, c)) &= W_{\rm cell} (F) 
\end{align*}
for all $R \in SO(d)$, $c \in \R^d$. 

\item $W_{\rm cell} (F)$ is minimal ($= 0$) if and only if there exists $R \in SO(d)$ and $c \in \R^d$ such that 
$$ F = R Z + (c, \ldots, c).$$ 

\item $W_{\rm cell}$ is $C^2$ in a neighborhood of $\bar{SO}(d) := SO(d) Z$ and the Hessian $D^2 W_{\rm cell}(Z)$ at the identity is positive definite on the orthogonal complement of the subspace spanned by translations $(c, \ldots, c)$ and infinitesimal rotations $FZ$, with $F^T = -F$. 

\item $p\geq d$, which together with Assumption \ref{ass:coercive} implies in particular that $W_{\rm cell}$ satisfies the growth assumption
\[ \liminf_{|F| \to \infty \atop F \in V_0} \frac{W_{\rm cell}(F)}{|F|^d} > 0. \]
\end{enumerate}
\end{ass}

\subsection{ $\Gamma$-convergence and integral representations}

In our analysis, we consider energies on discrete systems depending on a small parameter $\e$, the scale of the lattice spacing. To make the limit for $\e\to 0$ precise and gain some knowledge about the behavior of associated minimizers, we will use De Giorgi's $\Gamma$-convergence.
We recall the definition and some basic properties that will be needed in the sequel.
\begin{defi}
	Let $X$ be a metric space and $F_n,F \colon X \to \bar{\R}=\R\cup\{-\infty,\infty\} $. We say $F_n$ $\Gamma(X)$-converges to $F$ ($F_n \overset{\Gamma}{\longrightarrow} F$), if
	\begin{enumerate}[(i)]
		\item (liminf-inequality) For every $y,y_n\in X$ with $y_n\to y$, we have
		 \[ F(y) \leq \liminf\limits_{n\to\infty} F_n(y_n), \]
		\item (recovery sequence) For every $y\in X$, there is a sequence $y_n\in X$ such that
		 \[ F(y) \geq \limsup\limits_{n\to\infty} F_n(y_n). \] 
	\end{enumerate}
\end{defi}
If $(F_\e)_{\e>0}$ is a family of functionals depending on a positive real parameter $\e$, we say $F_\e$ $\Gamma(X)$-converges to $F$, if for every sequence $\e_n > 0$ converging to $0$, we have $F_{\e_n} \overset{\Gamma}{\longrightarrow} F$.
We will also use the $\Gamma$-$\limsup$ and the $\Gamma$-$\liminf$, given by
\al{
	F'(y) &= \Gamma(X)\text{-}\liminf\limits_{n\to\infty} F_n(y)= \inf\{ \liminf\limits_{n\to\infty} F_n(y_n) \colon y_n\to y \text{ in }  X \}, \\
	F''(y) &= \Gamma(X)\text{-}\limsup\limits_{n\to\infty} F_n(y)= \inf\{ \limsup\limits_{n\to\infty} F_n(y_n) \colon y_n\to y \text{ in }  X \}.
}
Note that (i) is equivalent to $F\leq F'$ and (ii) is equivalent to $F\geq F''$. Hence, $F_n \overset{\Gamma}{\longrightarrow} F$ if and only if $F'=F''=F$. Furthermore, we see that $\Gamma$-convergence is a pointwise property, so we can speak about $\Gamma$-convergence at a specific point.

In the following proposition we assemble some basic properties of $\Gamma$-convergence that we will not prove here.
\begin{prop} \label{prop:gammaproperties}
\begin{enumerate}[(i)]
	\item The infima in the definitions of $F'$ and $F''$ are actually attained minima in $\bar{\R}$;
	\item every sequence of functionals on a separable metric space, like $L^p(U;\R^d)$, has a $\Gamma$-convergent subsequence;
	\item $F'$, $F''$ and $F$ are lower semicontinuous with respect to convergence in $X$.
	\item $\Gamma$-convergence satisfies the Urysohn property, i.e., $F_n$ $\Gamma$-converges to $F$, if and only if every subsequence of $F_n$ has a further subsequence, that $\Gamma$-converges to $F$;
	\item if $F_n$ $\Gamma$-converges to $F$ and $G_n$ converges uniformly on bounded sets to a continuous functional $G$, then $F_n+G_n$ $\Gamma$-converges to $F+G$.
\end{enumerate}
\end{prop}
In view of applications, the most interesting property of $\Gamma$-convergence is the following theorem.
\begin{thm} \label{thm:minconverge}
	If $F_n$ $\Gamma$-converges to $F$ and sequences $y_n$ in $X$ with equibounded $F_n(y_n)$ are pre-compact 
	then $F$ attains its minimum on X and we have
	\[ \min\limits_{x\in X} F(x) = \lim\limits_{n \to \infty} \inf\limits_{x\in X} F_n(x). \]
	Furthermore, let $y_n\in X$ be a sequence with
	\[F_n(y_n) \to \lim\limits_{n \to \infty} \inf\limits_{x\in X} F_n(x),\]
	then the limit of every converging subsequence of $y_n$ is a minimizer of $F$.
\end{thm}
For proofs of Proposition \ref{prop:gammaproperties} and Theorem \ref{thm:minconverge} see, e.g., \cite{dalmasogamma}.

Returning to our specific setting, for a sequence $\e_n > 0$ such that $\e_n \to 0$, we define
\al{
	F'(y,U) &:= \Gamma(L^p(\Omega;\R^d))\text{-}\liminf\limits_{n\to\infty} F_{\e_n}(y,U) \\
	&= \min\{ \liminf\limits_{n\to\infty} F_{\e_n}(y_n,U) \colon y_n\to y \text{ in }  L^p(\Omega;\R^d)  \}, \\
	F''(y,U) &:= \Gamma(L^p(\Omega;\R^d))\text{-}\limsup\limits_{n\to\infty} F_{\e_n}(y,U). \\
	&= \min\{ \limsup\limits_{n\to\infty} F_{\e_n}(y_n,U) \colon y_n\to y \text{ in }  L^p(\Omega;\R^d)  \}.
}

The limiting functionals we will encounter in the next section are integral functionals of the form 
\[ I\colon W^{1,p}(U;\R^k)\to [0,\infty], \qquad  
I(y) = \int\limits_U f(\nabla y (x))\, dx \]
with $1<p<\infty$, $U\in\mathcal{A}(\R^d)$, $f\colon \R^{k\times d}\to [0,\infty)$ continuous. Recall that a Borel measurable and locally bounded function $f:\R^{k\times d}\to\R$ is quasiconvex, if
\[f(M) \leq \fint\limits_U f(M+\nabla \varphi (x)) \,dx,\]
for every nonempty $U\in\mathcal{A}(\R^d)$, $M\in\R^{k\times d}$ and $\varphi\in W^{1,\infty}_0(U;\R^k)$. In our analysis,  the quasiconvexity of $f$ will be due to the following result.
\begin{thm} \label{thm:quasiconvexlsc}
	If $I$ is sequentially weakly lower semicontinuous in $W^{1,p}(U;\R^k)$, then $f$ is quasiconvex.
\end{thm}

A detailed discussion of quasiconvexity and related properties, including proofs of the above statements, is given, e.g., in \cite{dacorogna}.

In order to guarantee that indeed our limiting functional is an integral functional, we will resort to the following general integral representation result on Sobolev spaces.
\begin{thm} \label{thm:integralrepresentation}
	Let $1\leq p < \infty$ and let $F\colon W^{1,p}(\Omega;\R^d)\times\mathcal{A}(\Omega) \to [0,\infty]$ satisfy the following conditions:
	\begin{enumerate}[(i)]
		\item (locality) $F(y,U) = F(v,U)$, if $y(x) = v(x)$ for a.e.\ $x\in U$;
		\item (measure property) $F(y,\cdot)$ is the restriction of a Borel measure to $\mathcal{A}(\Omega)$;
		\item (growth condition) there exists $c>0$ such that
			\[ F(y,U)\leq c \int\limits_U \abs{\nabla y(x)}^p + 1 \,dx;\]
		\item (translation invariance in y) $F(y,U) = F(y+a,U)$ for every $a\in\R^d$ ;
		\item (lower semicontinuity) $F(\cdot,U)$ is sequentially lower semicontinuous with respect to weak convergence in $W^{1,p}(\Omega;\R^d)$;
		\item (translation invariance in x) With $y_M(x)=M x$ we have
		\[ F(y_M,B_r(x))=F(y_M,B_r(x'))\]
		for every $M\in\R^{d\times d}$, $x,x'\in \Omega$ and $r>0$ such that $B_r(x),B_r(x')\subset\Omega$.
	\end{enumerate}
	Then there exists a continuous $f \colon \R^{d\times d} \to [0,\infty)$ such that
	\al{
		0\leq f(M) &\leq C(1+\abs{M}^p) \text{ for every } M\in\R^{d\times d}\text{ and} \\
		F(y,U) &= \int\limits_U f(\nabla y(x))\,dx .
	}
\end{thm}

A proof can be found in \cite[pp.77-81]{braidesdefranceschi} or in the scalar-valued setting, which is essentially the same, in \cite[pp.215-220]{dalmasogamma}.

To show the measure property in the previous theorem, we will use the following lemma. 
\begin{lem}[De Giorgi-Letta] \label{lem:degiorgiletta}
	Let $X$ be a metric space with open sets $\tau$. Assume that $\rho\colon \tau \to [0,\infty]$ is an increasing set function such that
	\begin{enumerate}[(i)]
		\item $\rho(\emptyset) = 0$,
		\item (subadditivity) $\rho(U\cup V) \leq \rho(U) + \rho (V)$ for all $U,V\in\tau$,
		\item (inner regularity) $\rho(U) = \sup\{\rho(V)\colon V\in\tau, V\css U\}$ for all $U\in\tau$,
		\item (superadditivity) $\rho(U\cup V) \geq \rho(U) + \rho (V)$ for all $U,V\in\tau$ with $U\cap V = \emptyset$.
	\end{enumerate}
	Then the extension $\mu$ of $\rho$ to all subsets of $X$, defined by
	\[
		\mu(E)= \inf\{ \rho(U)\colon U\in\tau,E\subset U \},  
	\]
	is an outer measure and every Borel set is $\mu$-measurable.	
\end{lem}
For a proof see, e.g., \cite[pp.32-34]{fonsecaleonilp}.

\section{A general representation result}\label{sec:repr}

In this section we will prove a general compactness and representation result for sequences of discrete deformations. For pair interactions, the following theorem has first been established by Alicandro and Cicalese in \cite{alicandrocicalese}.

\subsection{Statement of the representation result}

\begin{thm}[compactness and integral representation]
\label{thm:compactness}
Suppose Assumptions \ref{ass:coercive} and \ref{ass:growth} are true. For every sequence $\e_n > 0$ such that $\e_n \to 0$, there exists a subsequence $\e_{n_k}$ and a functional $F \colon L^p(\Omega;\R^d) \times \mathcal{A}(\Omega) \to [0,\infty]$ such that for every $U \in \mathcal{A}(\Omega)$ and $y\in W^{1,p}(\Omega;\R^d)$ the functionals $F_{\e_{n_k}}(\cdot,U)$ $\Gamma(L^p(\Omega;\R^d))$-converge to $F(\cdot,U)$ at $y$. Furthermore there exists a quasiconvex function $f\colon \R^{d\times d} \to [0,\infty)$ satisfying
\[c\abs{M}^p - c'\leq f(M) \leq c'(\abs{M}^p+1)\]
for some $c,c'>0$ such that
\[
	F(y,U)=\int\limits_{U} f(\nabla y(x))\,dx \quad \text{if } y\in W^{1,p}(\Omega;\R^d).
\]
In addition, if $U\in \mathcal{A}_L(\Omega)$ (in particular, if $U=\Omega$), we have
\[
	F(y, U)=
	\begin{dcases}
	 \int\limits_{U} f(\nabla y(x))\,dx & \text{if } y|_{U}\in W^{1,p}(U;\R^d),\\
	 \infty  & \text{otherwise,}
	\end{dcases}
\]
and the functionals $F_{\e_{n_k}}(\cdot,U)$ $\Gamma$-converge to $F(\cdot,U)$.
\end{thm}

\subsection{Interpolation}

We now define the continuous and piecewise affine interpolation $\tilde{y}$ of $y$, similar to \cite{schmidtlinelast}:
First consider the cell $A\left[-\frac{1}{2},\frac{1}{2}\right]^d$ and $y\colon A\left\{-\frac{1}{2},\frac{1}{2}\right\}^d \to \R^d$.
On every $0$-dimensional face of the cell just take $\tilde{y} = y$. Now assume we already have chosen a simplicial decomposition on every $(k-1)$-dimensional face and have interpolated affine there. Let $F = \co\{z_{i_1},\dots,z_{i_{2^k}}\}$ be a $k$-dimensional face.
Set
\[ \bar{z}=\frac{1}{2^k}\sum\limits_{m=1}^{2^k} z_{i_m}, \quad \tilde{y}(\bar{z}) =  \frac{1}{2^k}\sum\limits_{m=1}^{2^k} y(z_{i_m}) .\]
To complete the induction, we decompose $F$ into the simplices $\co \{w_1,\dots,w_k,\bar{z}\}$, where $\co \{w_1,\dots,w_k\}$ is a simplex belonging to a simplicial decomposition of an $(n-1)$-dimensional face. Define $\tilde{y}$ to be the interpolation affine on every constructed simplex. 
If $y \in \mathcal{B}_{\e}(\Omega)$, we get $\tilde{y}$ on $\overline{\Omega_\e}$ by interpolating as above on every cell.

The following proposition is about the relation of $\bar{\nabla}y$ and $\nabla{\tilde{y}}$.
\begin{prop} \label{prop:gradients}
There are $C,c > 0$ such that for every $x\in\Omega_\e$ and $y\in\mathcal{B}_{\e}(\Omega)$
\beq
	c \abs{\bar{\nabla}y(x)}^p  \leq \fint\limits_{Q_\e(x)} \abs{\nabla{\tilde{y}}(x')}^p \,dx' \leq C \abs{\bar{\nabla}y(x)}^p.
\eeq
\end{prop}
\bp
	Every $z_i$ belongs to some simplex $K$ of the construction. Choose $a\in K^\circ$, where the gradient is well-defined. Since the interpolation is linear on $K$, we see that
	\al{
		\abs{\frac{y(\bar{x}+\e z_i)-\bar{y}}{\e}}^p &=\abs{\nabla{\tilde{y}}(a)z_i}^p \\
		&\leq C \abs{\nabla{\tilde{y}}(a)}^p \\
		&= C\frac{1}{\abs{K}} \int\limits_K \abs{\nabla{\tilde{y}}(x')}^p\,dx' \\
		&\leq C \frac{\abs{Q_{\e}(x)}}{\abs{K}} \fint\limits_{Q_\e(x)} \abs{\nabla{\tilde{y}}(x')}^p\,dx' \\
		&\leq C \fint\limits_{Q_\e(x)} \abs{\nabla{\tilde{y}}(x')}^p\,dx',
	}
	where $C$ is independent of $x,\e$ and $y$.
	We immediately get the first inequality. For the second inequality we prove by induction over $k$ that for every $k$-dimensional simplex $S=\co\{\bar{z},z_i,w_1,\dots,w_{k-1}\}$ in the construction regarding $Q_\e(x)$ we have
	\beq \label{eq:induction1}
		\abs{\nabla{\tilde{y}}(a)P_{V}}^p\leq C \abs{\bar{\nabla}y(x)}^p
	\eeq
	for every $a\in S$, where $P_{V}$ is the projection on $V = \spano{ \{\bar{z}-z_i, w_1-z_i,\dots,w_{k-1}-z_i\} }$.
	The case $k=1$ is clear since then for some $j$ we have $V = \spano{ \{z_j-z_i\} }$ and \[\nabla{\tilde{y}}(a)(z_j-z_i) = \bar{\nabla}y(x)(e_j-e_i).\]
	If the statement is true for $k-1$, we immediately have \eqref{eq:induction1} for \[V' = \spano{ \{w_1-z_i,\dots,w_{k-1}-z_i\} }.\] But as in the $k=1$ case we also have \eqref{eq:induction1} for $V'' = \spano{\{\bar{z}-z_i\}}=\spano{\{z_j-z_i\}}$. Let us define
	\[ \norm{v}_V = \abs{v'}+ \abs{v''},\]
	if $v\in V$, $v'\in V'$ and $v''\in V''$ such that $v=v'+v''$.	This is a norm on $V$ and hence we can calculate using the equivalence of all norms on finite dimensional spaces
	\al{
		\abs{\nabla{\tilde{y}}(a)P_{V}}^p &\leq C \sup\{ \abs{\nabla{\tilde{y}}(a) v}^p \colon v\in V, \norm{v}_V\leq 1 \} \\
		&\leq C ( \sup\{ \abs{\nabla{\tilde{y}}(a) v'}^p \colon v'\in V', \abs{v'}\leq 1 \} + \sup\{ \abs{\nabla{\tilde{y}}(a) v'}^p \colon v'\in V', \abs{v'}\leq 1 \}) \\
		&\leq C \abs{\bar{\nabla}y(x)}^p.
	}
	Since we have only finite many possibilities for $V,V',V''$, this $C$ can be chosen independent of them, which concludes the induction.
	Take $k=d$ and integrate to get the result.
\ep

\begin{prop} \label{prop:interpol}
	Let $\e_n>0$, with $\e_n\to 0$, $y_n\in\mathcal{B}_{\e_n}(\Omega)$ and $y\in L^p(\Omega; \R^d)$ such that $y_n\to y$ in $L^p(\Omega; \R^d)$. For every $V\css\Omega$, we then have $\tilde{y}_n\to y$ in $L^p(V; \R^d)$.
\end{prop}
\bp
It is enough to show $\norm{y_n-\tilde{y}_n}_{L^p(V;\R^d)}\to 0$.
Let $\lambda_i \colon \R^d \to [0,1]$ denote the cell-periodic functions such that
	\al{ \tilde{y}_n (x) &= \sum\limits_{i=1}^{2^d} \lambda_i\left(\frac{x}{\e_n}\right) y_n(\bar{x} + \e_n z_i)\\
	&= \sum\limits_{i=1}^{2^d} \lambda_i\left(\frac{x}{\e_n}\right) y_n( x + \e_n (z_i-z_1)), }
	where without loss of generality we have chosen a numbering of $A \{-\frac{1}{2}, \frac{1}{2}\}^d$ such that $z_1 = A(-\frac{1}{2}, \ldots, -\frac{1}{2})$. Of course, $\lambda_i\geq 0$ and the $\lambda_i$ add up to $1$ in any point and so for $n$ large enough
\al{
	\int\limits_V \abs{y_n(x) - \tilde{y}_n(x)}^p\,dx &\leq \int\limits_V \left(\sum\limits_{i=1}^{2^d} \lambda_i\left(\frac{x}{\e_n}\right) \abs{y_n(x) - y_n(x+\e_n (z_i-z_1))}\right)^p\,dx \\
	&\leq \int\limits_V \max\limits_{i=1, \dotsc , 2^d}  \abs{y_n(x) - y_n(x+\e_n (z_i-z_1))}^p \,dx \\
	&\leq \sum\limits_{i=1}^{2^d} \int\limits_V  \abs{y_n(x) - y_n(x+\e_n (z_i-z_1))}^p \,dx.
}
But the last term goes to $0$ since for every $i\in\{1,\dotsc,2^d\}$
\al{
	\norm{y_n-y_n(\cdot + \e_n (z_i-z_1))}_{L^p(V;\R^d)} &\leq 2 \norm{y_n-y}_{L^p(\Omega;\R^d)} + \norm{y-y(\cdot + \e_n (z_i-z_1))}_{L^p(V;\R^d)} \\
	&\to 0.
}
\ep

\subsection{Preliminary lemmata}
We proceed to collect further lemmata. We will use them later to prove the requirements of Theorem \ref{thm:integralrepresentation}. In the following, fix some sequence of positive real numbers $\e_n \to 0$.

\begin{lem} \label{lem:gammaliminf}
 Suppose Assumption \ref{ass:coercive} is true. If $y \in L^p(\Omega;\R^d)$ and $U\in\mathcal{A}(\Omega)$ are such that $F'(y,U)<\infty$, then $y\in W^{1,p}(U;\R^d)$ and
 \beq
  F'(y,U) \geq c \norm{\nabla{y}}_{L^p(U;\R^{d \times d})}^p - c' \abs{U} ,
 \eeq
 for some $c,c' > 0$ independent of $y$ and $U$.
\end{lem}

\bp
	Let $y_n \to y$ in $L^p(\Omega;\R^d)$ such that $\liminf\limits_{n\to\infty} F_{\e_n}(y_n,U)<\infty$. For some subsequence $n_k$, we have
	\[\lim\limits_{k\to\infty} F_{\e_{n_k}}(y_{n_k},U) =\liminf\limits_{n\to\infty} F_{\e_n}(y_n,U), \]
	$y_{n_k}\in\mathcal{B}_{\e_{n_k}}(U)$ and $F_{\e_{n_k}}(y_{n_k},U) \leq M < \infty$ for some fixed $M>0$.
	By Proposition \ref{prop:interpol} we have $\tilde{y}_{n_k} \to y$ in $L^p(V, \R^d)$ for every $V \css U$. Furthermore, by Assumption \ref{ass:coercive} and Proposition \ref{prop:gradients}, we get
	\al{
		M \geq F_{\e_{n_k}}(y_{n_k},U) &= \e_{n_k}^d \sum\limits_{x\in (\la_{\e_{n_k}}'(U))^\circ} W_{\rm cell}(\bar{\nabla}y_{n_k}(x)) \\
		& \geq \e_{n_k}^d \sum\limits_{x\in (\la_{\e_{n_k}}'(U))^\circ} \left(c\abs{\bar{\nabla}y_{n_k}(x)}^p - c'\right) \\
		& \geq \e_{n_k}^d \sum\limits_{x\in (\la_{\e_{n_k}}'(U))^\circ} \left(c\fint\limits_{Q_{\e_{n_k}}(x)}\abs{\nabla\tilde{y}_{n_k}(x')}^p\,dx' - c'\right) .
	}
	We thus obtain
	\[c \int\limits_{U_{\e_{n_k}}}\abs{\nabla\tilde{y}_{n_k}(x')}^p\,dx' \leq M + c' \abs{U},\]
	hence the gradients are bounded in $L^p(V;\R^d)$. By the properties of weak convergence on Sobolev spaces this means $y\in W^{1,p}(V,\R^d)$ and $\nabla\tilde{y}_{n_k} \wto \nabla y$ in $L^p(V;\R^d)$. Weak sequentially lower semicontinuity of the norm yields
	\[c \norm{\nabla{y}}_{L^p(V;\R^{d \times d})}^p \leq \liminf\limits_{n\to\infty} F_{\e_n}(y_n,U) + c' \abs{U},\]
	but the right hand side is independent of $V$, thus $y\in W^{1,p}(U,\R^d)$ and
	\[c \norm{\nabla{y}}_{L^p(U;\R^{d \times d})}^p \leq \liminf\limits_{n\to\infty} F_{\e_n}(y_n,U) + c' \abs{U}.\]
	The definition of the $\Gamma$-$\liminf$ now yields the lemma.
\ep

\begin{lem} \label{lem:growth}
 Suppose Assumption \ref{ass:growth} is true. Then there is a $C>0$ such that for every $V\in\mathcal{A}_L(\Omega)$, $U\in\mathcal{A}(V)$ and $y\in L^p(\Omega;\R^d)\cap W^{1,p}(V;\R^d)$ we have
 \beq \label{eq:gammasupgrowth}
  F''(y,U) \leq C\left( \norm{\nabla{y}}_{L^p(U;\R^{d \times d})}^p + \abs{U} \right).
 \eeq
\end{lem}

\bp
	We first prove \eqref{eq:gammasupgrowth} for every $y\in C^{\infty}_c(\R^d;\R^d)$.
	For $x\in \la'_{\e_n}$ and $a\in Q_{\e_n}(x)$ define
	\[y_n(a)=y(x).\]
	Thus $y_n\in\mathcal{B}_{\e_n}(U)$ and since $y$ is uniformly continuous, we have $y_n \to y$ uniformly and hence in $L^p(\Omega;\R^d)$.
	By Taylor expansion we have
	\al{
		\abs{\frac{y_n(\bar{x}) - y_n(\bar{x}+\e_n z_i)}{\e_n}} &= \abs{\frac{y(\bar{x}) - y(\bar{x}+\e_n z_i)}{\e_n}} \\
		&\leq C ( \abs{\nabla{y}(\bar{x})} + \e_n \norm{\nabla^2 y}_{\infty} ) .
	}
	With Assumption \ref{ass:growth} we can calculate 
	\al{
		F_{\e_n}(y_n,U) &= \e_n^d \sum\limits_{x\in (\la_{\e_n}'(U))^\circ} W_{\rm cell}(\bar{\nabla}y_n(x)) \\
		&\leq C \e_n^d \sum\limits_{x\in (\la_{\e_n}'(U))^\circ} (\abs{\bar{\nabla}y_n(x)}^p + 1) \\
		&\leq C'\abs{U} + C\e_n^d \sum\limits_{x\in (\la_{\e_n}'(U))^\circ} ( \abs{\nabla{y}(x)}^p + \e_n^p \norm{\nabla^2 y}_{\infty}^p ) \\
		&\leq C'\abs{U} + C''\abs{U}\e_n^p \norm{\nabla^2 y}_{\infty}^p + C\e_n^d \sum\limits_{x\in (\la_{\e_n}'(U))^\circ} \abs{\nabla{y}(x)}^p.
	}	
	Furthermore for every $x'\in Q_{\e_n}(x)$, $x\in (\la_{\e_n}'(U))^\circ$
	\al{
		\abs{\nabla{y}(x)}^p &\leq C ( \abs{\nabla{y}(x')}^p + \abs{\nabla{y}(x)-\nabla{y}(x')}^p) \\
		&\leq C (\abs{\nabla{y}(x')}^p + \e_n^p \norm{\nabla^2 y}_\infty^p),
	}
	and, by integrating over $x'$ and summing over $x$, we get
	\al{ 
		\e_n^d \sum\limits_{x\in (\la_{\e_n}'(U))^\circ} \abs{\nabla{y}(x)}^p &\leq C \left( \int\limits_{U_{\e_n}} \abs{\nabla{y}(x')}^p\,dx' + \abs{U_{\e_n}} \e_n^p \norm{\nabla^2 y}_\infty^p \right) \\
		&\leq  C \left( \int\limits_U \abs{\nabla{y}(x')}^p\,dx' + \abs{U} \e_n^p \norm{\nabla^2 y}_\infty^p \right).
	}
	Putting the two inequalities together and letting $n\to\infty$, we obtain
	\[ \limsup\limits_{n\to\infty} F_{\e_n}(y_n,U) \leq  C\left( \norm{\nabla{y}}_{L^p(U;\R^{d \times d})}^p + \abs{U} \right). \]
	So by the definition of the $\Gamma$-$\limsup$ we have \eqref{eq:gammasupgrowth}.
	
	The general case follows easily: Since $V$ has Lipschitz boundary, we can take $y_k\in C^{\infty}_c(\R^d;\R^d)$ such that $y_k \to y$ in $W^{1,p}(V;\R^d)$. Then we have by lower semicontinuity of $F''(\cdot,U)$
	\al{
		F''(y,U) &\leq \liminf\limits_{k\to\infty} F''(y_k,U) \\
		&\leq \liminf\limits_{k\to\infty} C\left( \norm{\nabla{y_k}}_{L^p(U;\R^{d \times d})}^p + \abs{U} \right) \\
		& = C\left( \norm{\nabla{y}}_{L^p(U;\R^{d \times d})}^p + \abs{U} \right).
	}
\ep

\begin{lem} \label{lem:subadd}
	Suppose Assumptions \ref{ass:coercive} and \ref{ass:growth} are true. Let $U,V,U' \in\mathcal{A}(\Omega)$ be such that $U'\css U$. Then for every $y\in W^{1,p}(\Omega;\R^d)$
	\[ F''(y,U'\cup V) \leq F''(y,U) + F''(y,V). \]
\end{lem}
\bp
	Without loss of generality, we can assume the terms on the right hand side to be finite. According to the properties of the $\Gamma$-$\limsup$ it is possible to find sequences $u_n,v_n$ such that
	\al{
		&\limsup\limits_{n\to\infty} F_{\e_n}(u_n,U) = F''(y,U)\\
		&\limsup\limits_{n\to\infty} F_{\e_n}(v_n,V) = F''(y,V)\\
		&u_n\to y \text{ in } L^p(\Omega;\R^d)\\
		&v_n\to y \text{ in } L^p(\Omega;\R^d)
	}
	For $n$ large enough $F_{\e_n}(u_n,U)$ and $F_{\e_n}(v_n,V)$ are bounded and $u_n\in\mathcal{B}_{\e_n}(U)$, $v_n\in\mathcal{B}_{\e_n}(V)$.
	
	Fix $N\in\N$, $N\geq 5$ and then define $D=\dist(U',U^c)$ and $U_j = \{x\in U \colon \dist(x,U')<\frac{j D}{N}\}$. Choose cut-off functions $\varphi_j$ such that
	\al{
		&\varphi_j(x) = 1\quad \forall x\in U_j, \\
		&\varphi_j\in C_c^\infty(U_{j+1};[0,1]), \\
		&\norm{\nabla\varphi_j}_\infty \leq \frac{2 N}{D}.
	}
	Next we define
	\[ w_{n,j}(x) = \varphi_j(\bar{x}) u_n(x) + (1 - \varphi_j(\bar{x})) v_n(x) \]
	and calculate
	\begin{align} \label{eq:diffquoeq}
		\frac{w_{n,j}(x+\e_n z_i) - w_{n,j}(x)}{\e_n} =\ &\varphi_j(\overline{x+\e_n z_i}) \frac{u_n(x+\e_n z_i) - u_n(x)}{\e_n} \nonumber \\
		&+ (1 - \varphi_j(\overline{x+\e_n z_i})) \frac{v_n(x+\e_n z_i) - v_n(x)}{\e_n} \\
		&+ (u_n(x) - v_n(x)) \frac{\varphi_j(\overline{x+\e_n z_i}) - \varphi_j(\bar{x})}{\e_n} \nonumber .
	\end{align}
	To estimate $F_{\e_n}(w_{n,j},U'\cup V)$, we have to look at $(\la_{\e_n}'(U'\cup V))^\circ$. Clearly, if $x$ is in $(\la_{\e_n}'(U_j))^\circ$, then $\bar{\nabla}w_{n,j}(x) = \bar{\nabla}u_n(x) $ and if $x$ is in $(\la_{\e_n}'(V\backslash\overline{U_{j+1}}))^\circ$, then $\bar{\nabla}w_{n,j}(x) = \bar{\nabla}v_n(x)$. To control the other cases, observe that for $n$ large enough $\diam(Q_{\e_n})\leq \frac{D}{2 N}$ and thus
	\[ (\la_{\e_n}'(U'\cup V))^\circ \subset (\la_{\e_n}'(U_j))^\circ \cup (\la_{\e_n}'(V\backslash\overline{U_{j+1}}))^\circ \cup (\la_{\e_n}'(V\cap(U_{j+2}\backslash\overline{U_{j-1}})))^\circ  \]
	for every $j\in\{2,\dots,N-3\}$ and $n$ large enough.
	With $W_j = V\cap(U_{j+2}\backslash\overline{U_{j-1}})$, we then have
	\al{
		F_{\e_n}(w_{n,j},U'\cup V) &= \e_n^d \sum\limits_{x\in (\la_{\e_n}'(U'\cup V))^\circ} W_{\rm cell}(\bar{\nabla}w_{n,j}(x)) \\
		&\leq F_{\e_n}(u_n,U) + F_{\e_n}(v_n,V) + \underbrace{ \e_n^d \sum\limits_{x\in (\la_{\e_n}'(W_j))^\circ} W_{\rm cell}(\bar{\nabla}w_{n,j}(x)) }_{:= S_{j,n}} .
	}
	We now have to estimate $S_{j,n}$. For all $n$ large enough, use first Assumption \ref{ass:growth} and then \eqref{eq:diffquoeq} to get
	\al{
		S_{j,n} &\leq C \e_n^d \sum\limits_{x\in (\la_{\e_n}'(W_j))^\circ} ( \abs{\bar{\nabla}w_{n,j}(x)}^p + 1) \\
		&\leq C \e_n^d \sum\limits_{x\in (\la_{\e_n}'(W_j))^\circ} ( \abs{\bar{\nabla}u_n(x)}^p + \abs{\bar{\nabla}v_n(x)}^p + \abs{u_n(x) - v_n(x)}^p \norm{\nabla\varphi_j}_\infty^p + 1) \\
		&\leq C \e_n^d \sum\limits_{x\in (\la_{\e_n}'(W_j))^\circ} ( \abs{\bar{\nabla}u_n(x)}^p + \abs{\bar{\nabla}v_n(x)}^p + \abs{u_n(x) - v_n(x)}^p N^p + 1) \\
		&\leq C \int\limits_{(W_j)_{\e_n}}  \abs{\nabla\tilde{u}_n(x)}^p + \abs{\nabla\tilde{v}_n(x)}^p + N^p \abs{u_n(x) - v_n(x)}^p + 1 \,dx ,
	}
	because of the gradient of $\varphi$ being bounded by $C N$ and Proposition \ref{prop:gradients}.
	Averaging over $j$, we get
	\beq \label{eq:avS}
		\frac{1}{N-4}\sum\limits_{j=2}^{N-3} S_{j,n} \leq C \frac{1}{N-4} \int\limits_{V_{\e_n}}  \abs{\nabla\tilde{u}_n(x)}^p + \abs{\nabla\tilde{v}_n(x)}^p + 1 \,dx + N^p \int\limits_{V_{\e_n}} \abs{u_n(x) - v_n(x)}^p \,dx .
	\eeq
Of course we can always find a number $j(n)$ such that
	\[S_{j(n),n}\leq\frac{1}{N-4}\sum\limits_{j=2}^{N-3} S_{j,n}.\]
	By Proposition \ref{prop:gradients} and Assumption \ref{ass:coercive}, the first integral in \eqref{eq:avS} is bounded, but 
	\[\norm{u_n-v_n}_{L^p(\Omega;\R^d)}\to 0\] for $n\to\infty$, hence
	\[ \limsup\limits_{n\to\infty} S_{j(n),n} \leq \frac{C}{N-4}. \]
	If we define $y_n = w_{n,j(n)}$, then obviously $y_n\in\mathcal{B}_{\e_{n}}(U'\cup V)$ and $y_n \to y$ in $L^p(\Omega;\R^d)$. We have
	\al{
		F''(y,U'\cup V) &\leq \limsup\limits_{n\to\infty} F_{\e_n}(y_n,U'\cup V) \\
		&\leq \limsup\limits_{n\to\infty} F_{\e_n}(u_n,U) + \limsup\limits_{n\to\infty} F_{\e_n}(v_n,V) +\limsup\limits_{n\to\infty} S_{j(n),n} \\
		&\leq	F''(y,U) + F''(y,V) + \frac{C}{N-4}.
	}
	Letting $N\to\infty$, we get the conclusion. 
\ep

\begin{lem} \label{lem:limitregular}
	Suppose Assumptions \ref{ass:coercive} and \ref{ass:growth} are true. Then for every $V\in\mathcal{A}_L(\Omega)$, $U\in\mathcal{A}(V)$ and $y\in L^p(\Omega;\R^d)\cap W^{1,p}(V;\R^d)$
	\[ F''(y,U) = \sup\limits_{U'\css U} F''(y,U'). \]
\end{lem}

\bp
Since $F''(y,\cdot)$ is an increasing set function, we only have to show '$\leq$'.

Let $\delta>0$. Then take a $U'''\css U$ such that
\[ \abs{U\backslash \overline{U'''}} +\norm{\nabla y}_{L^p(U\backslash \overline{U'''};\R^d)}\leq \delta .\]
Choosing $U',U''$ such that
\[ U''' \css U'' \css U' \css U, \]
we can calculate
\al{ F''(y,U) &\leq F''(y, U'' \cup U\backslash \overline{U'''} ) \\
&\leq F''(y,U') + F''(y, U\backslash \overline{U'''} ) \\
&\leq F''(y,U') + \delta C,
}
where we used Lemma \ref{lem:subadd} and Lemma \ref{lem:growth}.
\ep

\begin{lem} \label{lem:locality}
	Suppose Assumptions \ref{ass:coercive} and \ref{ass:growth} are true. Then for every $V\in\mathcal{A}_L(\Omega)$, $U\in\mathcal{A}(V)$ and $u,v\in L^p(\Omega;\R^d)\cap W^{1,p}(V;\R^d)$ such that $u(x) = v(x)$ for almost every $x\in U$, we have
	\[ F''(u,U) = F''(v,U).\]
\end{lem}

\bp
If $u=v$ a.e.\ in $U$ then for $U'\css U$ we have $F''(u,U') = F''(v,U')$. To see this, just change any approximating discrete sequence of $u$ outside of $(U')^{\e_n}$ such that the new sequence converges to $v$. But this is enough by Lemma \ref{lem:limitregular}.
\ep

\subsection{Proof of the representation result}

Now, we can finally prove the compactness result:

\bp[Proof of Theorem \ref{thm:compactness}]
First we find by a suitable diagonal argument a subsequence $F_{\e_{n_k}}$ such that we get $\Gamma$-convergence for every $U\in\mathcal{A}(\Omega)$. For this we define
\[ \mathcal{A}_1 = \left\{ U\subset\Omega \colon U = \bigcup\limits_{i=1}^N B_{r_i}(x_i), x_i\in\Q^d, r_i\in\Q, r_i > 0, N\in\N . \right\} \] 
The set $\mathcal{A}_1$ is countable and we can write $\mathcal{A}_1 =\{U_1,U_2,\dots\}$. Now choose subsequences as follows:
\[
\begin{array}{ccc}
	F_{\e_n}(\cdot,U_1) &\text{ has a } \Gamma\text{-convergent subsequence } &F_{\e_{n^1_k}}(\cdot,U_1), \\
	F_{\e_{n^1_k}}(\cdot,U_2) &\text{ has a } \Gamma\text{-convergent subsequence } &F_{\e_{n^2_k}}(\cdot,U_2), \\
	F_{\e_{n^2_k}}(\cdot,U_3) &\text{ has a } \Gamma\text{-convergent subsequence } &F_{\e_{n^3_k}}(\cdot,U_3), \\
	\vdots & \vdots &\vdots \\
\end{array}
\]
Now setting $n_k = n_k^k$, we see that $F_{\e_{n_k}}(\cdot,U)$ $\Gamma$-converges to a $F(\cdot,U)$ for every $U\in\mathcal{A}_1$. In the following we will only consider the sequence $\e_{n_k}$ and, in particular, define $F'$ and $F''$ accordingly. Furthermore, we define $F(y,U):=F'(y,U)$ for every $y$ and $U$.

For $W\css U \subset \Omega$, by compactness of $\overline{W}$, we always find $V\in\mathcal{A}_1$ such that $W\subset V \css U$. Hence, by Lemma \ref{lem:limitregular} we have
\[ F''(y,U) = \sup\{ F''(y,V) \colon V\css U, V\in\mathcal{A}_1 \} \]
for every $U\in\mathcal{A}(\Omega)$ and $y\in W^{1,p}(\Omega,\R^d)$. Using, that $F'(y,\cdot)$ is an increasing set function, we can calculate
\al{
	\sup\{ F'(y,V) \colon V\css U, V\in\mathcal{A}_1 \} &\leq F'(y,U) \\
	&\leq F''(y,U) \\
	&= \sup\{ F''(y,V) \colon V\css U, V\in\mathcal{A}_1 \}.
}
But the first and the last term are equal, thus $F'(y,U) = F''(y,U) = F(y,U)$, whenever $y\in W^{1,p}(\Omega,\R^d)$.

The next step is to get an integral representation by showing that $F$, restricted to $W^{1,p}(\Omega;\R^d)$, satisfies the conditions (i)-(vi) in Theorem \ref{thm:integralrepresentation}. We immediately see the locality (i), by Lemma \ref{lem:locality}, and the growth condition (iii), by Lemma \ref{lem:growth}. Furthermore, since the $F_{\e_{n_k}}$ are translation invariant in $y$, so is $F$, which yields (iv). To get the lower semicontinuity (v), just remember that weak convergence in $W^{1,p}(\Omega,\R^d)$ implies strong convergence in $L^p(\Omega,\R^d)$ and that $\Gamma$(X)-limits are sequentially lower semicontinuous with respect to the convergence in $X$.

To get the measure property (ii), it is enough to show that we can apply the De-Giorgi-Letta criterion (Lemma \ref{lem:degiorgiletta}) with $\rho = F(y,\cdot)$.
Obviously $F(y,\cdot)$ is an increasing set function and $F(y,\emptyset)=0$. Remark that for every $W\css U \cup V$($W,U,V$ open), there are open sets $U',V'$ such that $U'\css U$, $V'\css V$ and $W\subset U'\cup V'$, which is easily seen by the compactness of $\overline{W}$. Hence the subadditivity follows from the Lemmata \ref{lem:subadd} and \ref{lem:limitregular}. The inner regularity is explicitly given by Lemma \ref{lem:limitregular}. The superadditivity we can show directly. Take a sequence $y_k\in\mathcal{B}_{\e_{n_k}}(U\cup V)$ such that $y_k\to y$ in $L^p(\Omega;\R^d)$ and
\[ F(y,U\cup V) = \lim\limits_{k\to\infty} F_{\e_{n_k}}(y_k,U \cup V).\]
Then,
\al{
	F(y,U\cup V) &\geq \liminf\limits_{k\to\infty} F_{\e_{n_k}}(y_k,U) + \liminf\limits_{k\to\infty} F_{\e_{n_k}}(y_k,V) \\
	&\geq F(y,U) + F(y,V),
}
since $U\cap V = \emptyset$. Hence, we can apply the De-Giorgi-Letta criterion and obtain (ii).
Finally, condition (vi) states that for every $M\in\R^{d\times d}$, $z,z'\in \Omega$ and $r>0$ such that $B_r(z),B_r(z')\subset\Omega$, we have
\[ F(y_M,B_r(z))=F(y_M,B_r(z')),\]
if we set $y_M(x) = Mx$. By inner regularity, it is enough to show that, for any $r' < r$,
\[F(y_M,B_r(z))\geq F(y_M,B_{r'}(z')).\]
Let $y_k\in\mathcal{B}_{\e_{n_k}}(B_r(z))$ such that $y_k \to y_M$ in $L^p(\Omega;\R^d)$ and
\[ \lim\limits_{k\to\infty} F_{\e_{n_k}}(y_k,B_r(z)) = F(y_M,B_r(z)).\]
Denote by $a_k$ the only point in $\la_{\e_{n_k}}\cap Q_{\e_{n_k}}(z'-z)$. Then define

\[
 u_k(x) = 
 \begin{dcases}
  y_k(x-a_k) + M a_k & \text{if } x\in (B_{r'}(z'))^{\e_{n_k}}\\
  M \bar{x}  & \text{else.}
 \end{dcases}
\]
If $k$ is large enough, then $x-a_k \in (B_r(z))_{\e_{n_k}}$, $u_k\in\mathcal{B}_{\e_{n_k}}(B_{r'}(z'))$ and 
\[ \bar{\nabla} u_k(x) = \bar{\nabla} y_k(x - a_k)\]
for all $x\in (B_{r'}(z'))^{\e_{n_k}}$. Hence,
\[F_{\e_{n_k}}(u_k,B_{r'}(z')) \leq F_{\e_{n_k}}(y_k,B_r(z)). \]
Furthermore, we have $M a_k \to M(z'-z)$ and $y_k(\cdot-a_k)\to M(\cdot-(z'-z))$ in $L^p(B_{r'}(z');\R^d)$ and therefore $u_k\to y_M$ in $L^p(\Omega;\R^d)$.
Hence, we get
\al{
	F(y_M,B_{r'}(z')) &\leq \liminf\limits_{k\to\infty} F_{\e_{n_k}}(u_k,B_{r'}(z')) \\
	&\leq  \liminf\limits_{k\to\infty}F_{\e_{n_k}}(y_k,B_r(z)) \\
	& = F(y_M, B_r(z)),
}
and (vi) is proven.

Consequently, we can apply Theorem \ref{thm:integralrepresentation} to the restriction of $F$ to $W^{1,p}(\Omega;\R^d)\times\mathcal{A}(\Omega)$.
In particular, there is a continuous function $f \colon \R^{d\times d} \to [0,\infty)$ such that
\[
 F(y,U) = \int\limits_U f(\nabla y(x))\,dx \quad \text{if } y\in W^{1,p}(\Omega;\R^d)  
\]
and
\beq
	0\leq f(M) \leq C(1+\abs{M}^p) \text{ for every } M\in\R^{d\times d}.
\eeq
The asserted lower bound on $f$ is instantly obtained, if we apply Lemma \ref{lem:gammaliminf} to $y_M$ and use the integral representation.
And finally, $f$ is quasiconvex by Theorem \ref{thm:quasiconvexlsc}, since $F(\cdot,\Omega)$ is sequentially lower semicontinuous with respect to weak convergence in $W^{1,p}(\Omega;\R^d)$.

Now, let $U$ have Lipschitz boundary. Take $y\in L^p(\Omega;\R^d)\cap W^{1,p}(U,\R^d)$. By Lemma \ref{lem:limitregular}, we have
\[ F''(y,U) = \sup\{ F''(y,V) \colon V\css U, V\in\mathcal{A}_1 \}. \]
Using that $F'(y,\cdot)$ is an increasing set function, we can calculate
\al{
	\sup\{ F'(y,V) \colon V\css U, V\in\mathcal{A}_1 \} &\leq F'(y,U) \\
	&\leq F''(y,U)\\
	&= \sup\{ F''(y,V) \colon V\css U, V\in\mathcal{A}_1 \}.
}
But the first and the last term are equal, thus $F'(y,U) = F''(y,U)= F(y,U)$.
If $y\in L^p(\Omega;\R^d)\backslash W^{1,p}(U,\R^d)$, then $\infty = F'(y,U) = F''(y,U) = F(y,U)$ by Lemma \ref{lem:gammaliminf}. Hence, $F_{\e_{n_k}}(\cdot,U)$ $\Gamma(L^p(\Omega;\R^d))$-converges to $F(\cdot,U)$. 
To get the integral representation for $y\in L^p(\Omega;\R^d)\cap W^{1,p}(U,\R^d)$, observe, that since $U$ has Lipschitz boundary, we can find a function $v\in W^{1,p}(\Omega;\R^d)$ such that
\[y(x) = v(x) \quad \text{for almost every } x\in U.\]
Then, by Lemma \ref{lem:locality},
\[ F(y,U) = F(v,U) = \int\limits_{U} f(\nabla v(x))\,dx = \int\limits_{U} f(\nabla y(x))\,dx .\] 
\ep

\section{The boundary value problem} \label{sec:bvp}

While loading terms can be included in our results so far without difficulties, the restriction to deformations with preassigned boundary values is more subtle. 

\subsection{Statement of representation result with boundary conditions}

Suppose $g\in W^{1,\infty}(\R^d;\R^d)$ is a boundary datum. We will then always choose the precise representative for $g$ and thus assume that $g$ is continuous. We define the admissible lattice deformations $\mathcal{B}_\e(U,g)$ as the functions in $\mathcal{B}_{\e}(U)$, that satisfy the boundary condition
\[ y(x) = g(\bar{x}) \text{, whenever } x\in \partial\la_\e(U).\]
The correspondingly restricted discrete functional is 
\[
	F_\e^g(y,U)=
	\begin{dcases}
		F_\e(y,U) & \text{ if } y\in\mathcal{B}_\e(U,g), \\
		\infty & \text{ otherwise.}
	\end{dcases}
\]
Assume that $\e_{n_k}$ and $f$ are as in Theorem \ref{thm:compactness}, let us for simplicity write just $\e_k$ in the following and set 
\[
	F^g(y,U)=
	\begin{dcases}
		F(y,U) & \text{ if } y|_U\in g + W_0^{1,p}(U;\R^d), \\
		\infty & \text{ otherwise.}
	\end{dcases}
\]
In analogy to Theorem \ref{thm:compactness} we then have: 
\begin{thm} \label{thm:boundarygamma}
Suppose Assumptions \ref{ass:coercive} and \ref{ass:growth} are true, $g\in W^{1,\infty}(\R^d;\R^d)$ and $F^g$, $F^g_{\e_k}$ are as above. Then $F^g_{\e_k}(\cdot, U)$ $\Gamma(L^p(\Omega;\R^d))$-converges to $F^g(\cdot,U)$ for every $U\in\mathcal{A}_L(\Omega)$.
\end{thm}

\subsection{Improved estimates on interpolations}

We start by improving Proposition \ref{prop:interpol} for sequences in $\mathcal{B}_\e(U,g)$. This is possible, because now we can control what happens near the boundary. Note, that now we can naturally define the interpolation $\tilde{y}$ on all of $U$, namely, we just extend $y$ by the discretization of $g$ before we interpolate. 

\begin{prop} \label{prop:interpolbc}
	Let $U\in\mathcal{A}_L(\Omega)$ and $y_k \in \mathcal{B}_{\e_k}(U,g)$. Then $y_k\to y$ in $L^p(U;\R^d)$ if and only if $\tilde{y}_k\to y$ in $L^p(U;\R^d)$.
\end{prop}
\bp
	First, let $y_k\to y$ in $L^p(U;\R^d)$. Choose some open bounded set $U'$ with Lipschitz boundary and $U\css U'$. Extend the functions by defining $y_k(x) := g(\bar{x})$ and $y(x) := g(x)$ for $x\in U'\backslash U$. So $y_k\in \mathcal{B}_{\e_k}(U',g)$ and, since $g$ is Lipschitz, we have $y_k\to y$ in $L^p(U';\R^d)$. But then by Proposition \ref{prop:interpol} we get $\tilde{y}_k\to y$ in $L^p(U;\R^d)$. 
	
Now, let $\tilde{y}_k\to y$ in $L^p(U;\R^d)$. Let $\lambda_i \colon \R^d \to [0,1]$ again denote the cell-periodic functions such that, with $z_1 = A(-\frac{1}{2}, \ldots, -\frac{1}{2})$, 
	\[ \tilde{y}_n (x) = \sum\limits_{i=1}^{2^d} \lambda_i\left(\frac{x}{\e_n}\right) y_n( x + \e_n (z_i-z_1)). \]
	Of course, $\lambda_i\geq 0$ and the $\lambda_i$ add up to $1$ in any point. Define
	\[W_{n,i} = \left\{x\in U \colon \lambda_i\left(\frac{x}{\e_n}\right)\geq \frac{1}{2} \text{ and } \lambda_j\left(\frac{x}{\e_n}\right)\leq a \text{ for } j\neq i \right\}, \]
	where $a$ will be chosen suitably later, and note that, for every $x \in U_{\e_n}$, the ratio $\frac{|W_{n,i} \cap Q_{\e_n}(x)|}{|Q_{\e_n}(x)|}$ is independent of $n$ and $x$ and positive since $x \in Q_{\e_n}$, $x \to \bar{x} + z_i$ implies that $\lambda_i(\frac{x}{\e_n}) \to 1$ and $\lambda_j(\frac{x}{\e_n}) \to 0$ for $j \neq i$. Next, extend $y$ and $y_k$ by $g$ as above and define
\[ P_n y (x):= \fint\limits_{Q_{\e_n}(x)} y(b)\,db.  \]
Of course, we have $\norm{P_n y - y}_{L^p(U;\R^d)}\to 0$. Hence, it suffices to show $\norm{P_n y - y_n}_{L^p(U;\R^d)}\to 0$.
For $x\in W_{n,i}$ we have
\al{
	\abs{\tilde{y}_n( x-\e_n (z_i-z_1)) - P_n y(x)} \geq &\frac{1}{2} \abs{y_n(x) - P_n y(x)} \\ &- \sum\limits_{j\neq i} \lambda_j\left(\frac{x}{\e_n}\right) \abs{y_n(x - \e_n (z_i-z_j)) - P_n y(x)}.
}
Since $y_n$ and $P_n y$ are constant on every cell, we thus have 
\al{
	\frac{1}{2}\norm{y_n-P_n y}_{L^p(U_{\e_n};\R^d)} &=  \frac{1}{2} \frac{\abs{U_{\e_n}}^{\frac{1}{p}}}{\abs{U_{\e_n}\cap W_{n,i}}^{\frac{1}{p}}}\norm{y_n-P_n y}_{L^p(U_{\e_n}\cap W_{n,i};\R^d)} \\
	&\leq \frac{\abs{U_{\e_n}}^{\frac{1}{p}}}{\abs{U_{\e_n}\cap W_{n,i}}^{\frac{1}{p}}} \norm{\tilde{y}_n(\cdot-\e_n (z_i-z_1)) - P_n y}_{L^p(U_{\e_n}\cap W_{n,i};\R^d)} \\
	&\quad + a \sum\limits_{j\neq i} \norm{y_n(\cdot-\e_n (z_i-z_j)) - P_n y}_{L^p(U_{\e_n};\R^d)}.
}
But $\frac{\abs{U_{\e_n}}^{\frac{1}{p}}}{\abs{U_{\e_n}\cap W_{n,i}}^{\frac{1}{p}}} > 0$ is independent of $n$ and  
\al{
	&\norm{\tilde{y}_n(\cdot-\e_n (z_i-z_1)) - P_n y}_{L^p(U_{\e_n}\cap W_{n,i};\R^d)} \\
	&\leq \norm{y - P_n y}_{L^p(U;\R^d)} + \norm{\tilde{y}_n(\cdot-\e_n (z_i-z_1)) - y}_{L^p(U;\R^d)} 
}
converges to $0$. To control the remaining sum, we estimate
\al{
	&\norm{y_n(\cdot-\e_n (z_i-z_j)) - P_n y}_{L^p(U_{\e_n};\R^d)} \\
	&\leq \norm{y_n(\cdot-\e_n (z_i-z_j)) - P_n y(\cdot-\e_n (z_i-z_j))}_{L^p(U_{\e_n};\R^d)} \\
	&\quad + \norm{P_n y(\cdot-\e_n (z_i-z_j)) - P_n y}_{L^p(U_{\e_n};\R^d)},
}
where the second term goes to $0$ and the first term is estimated by 
\al{ 
	\norm{y_n(\cdot-\e_n (z_i-z_j)) - P_n y(\cdot-\e_n (z_i-z_j))}_{L^p(U_{\e_n};\R^d)} 
	\leq \norm{y_n-P_n y}_{L^p(U^{\e_n};\R^d)}.
}
Altogether we obtain
\al{
	&\norm{y_n-P_n y}_{L^p(U;\R^d)} \\
	&\leq \norm{y_n-P_n y}_{L^p(U^{\e_n}\backslash U_{\e_n};\R^d)} + \norm{y_n-P_n y}_{L^p(U_{\e_n};\R^d)} \\ 
	&\leq \norm{y_n-P_n y}_{L^p(U^{\e_n}\backslash U_{\e_n};\R^d)} 
	+ 2a(2^d-1) \norm{y_n-P_n y}_{L^p(U^{\e_n};\R^d)} + o(1) \\ 
	&= (1+2a(2^d-1))\norm{y_n-P_n y}_{L^p(U^{\e_n}\backslash U_{\e_n};\R^d)} 
	+ 2a(2^d-1) \norm{y_n-P_n y}_{L^p(U_{\e_n};\R^d)} + o(1)
}
As near the boundary we can calculate
\al{
	\norm{y_n-P_n y}_{L^p(U^{\e_n}\backslash U_{\e_n};\R^d)} 
	&\leq \abs{U^{\e_n}\backslash U_{\e_n}}^{\frac{1}{p}}\norm{g}_{\infty} + \norm{P_n y}_{L^p(U^{\e_n}\backslash U_{\e_n};\R^d)} \\
	&\leq \abs{U^{\e_n}\backslash U_{\e_n}}^{\frac{1}{p}}\norm{g}_{\infty} + \norm{y}_{L^p(U^{\e_n}\backslash U_{\e_n};\R^d)} \\
	&\to 0, 
}
for $a = \frac{1}{2^{d+1}}$ we finally get
\[\norm{y_n-P_n y}_{L^p(U;\R^d)}\to 0.\]
\ep

\begin{rem}\label{rem:Lploc-conv}
The proof shows that without boundary conditions, i.e., for a general sequence $y_k \in \mathcal{B}_{\e_k}(U)$, $U\in\mathcal{A}(\Omega)$, we still have $y_k\to y$ in $L^p_{loc}(U;\R^d)$ if and only if $\tilde{y}_k\to y$ in $L^p_{loc}(U;\R^d)$. 
\end{rem}

\subsection{Proof of the boundary value representation result}

\bp[Proof of Theorem \ref{thm:boundarygamma}]
	Fix $U\in \mathcal{A}_L(\Omega)$.
	We start with the $\liminf$-inequality. Let $y_k,y\in L^p(\Omega;\R^d)$ such that $y_k \to y$. We can assume that
	\[ \liminf\limits_{k\to\infty} F^g_{\e_k}(y_k,U) < \infty,\]
	because otherwise there is nothing to show. For some subsequence we then get
	\al{
		\liminf\limits_{k\to\infty} F^g_{\e_k}(y_k,U) = \lim\limits_{l\to\infty} F^g_{\e_{k_l}}(y_{k_l},U). 		
	}
	But since $F_{\e_{k_l}} \leq F^g_{\e_{k_l}}$, we can argue as in Lemma \ref{lem:gammaliminf} to see that $y\in W^{1,p}(U;\R^d)$ and, for any $V\css U$, that $\tilde{y}_{k_l}\wto y$ in $W^{1,p}(V;\R^d)$. Using Proposition \ref{prop:interpolbc}, we see that $\tilde{y}_{k_l}$ converges strongly in $L^p(U;\R^d)$ and, since $\nabla\tilde{y}_{k_l}$ is now bounded in $L^p(U;\R^d)$, weakly in $W^{1,p}(U;\R^d)$ to $y$. Regarding the boundary condition, there are open neighborhoods $V_l$ of $\partial U$, where $\tilde{y}_{k_l}$ is an affine interpolation of $g$. Namely, $V_l$ is the interior of the union of all cells $Q_{\e_{k_l}}$, with $\overline{Q_{\e_{k_l}}}\cap \partial U \neq \emptyset$. Then 
\[\sup_{x \in \partial U}|\tilde{y}_{k_l}(x) - g(x)| 
\le \sup_{x \in V_l}|\tilde{y}_{k_l}(x) - g(x)|
\le C \e_{k_l} \]
since $g$ is Lipschitz. Denoting the trace operator by $T$, we thus have $T\tilde{y}_{k_l}\wto Ty = Tg$ in $L^p(\partial U;\R^d)$ and hence $y\in g + W_0^{1,p}(U;\R^d)$.
	But then, we can calculate
	\[F^g(y,U) = F(y,U)\leq \liminf\limits_{k\to\infty} F_{\e_k}(y_k,U) \leq \liminf\limits_{k\to\infty} F^g_{\e_k}(y_k,U), \]
	and have indeed proven the $\liminf$-inequality.

	To get the $\Gamma$-convergence result, we now proof the $\limsup$-inequality.
	Let us first assume $y(x) = g(x) + \psi(x)$, for every $x\in U$ and some $\psi \in C_c^\infty(U;\R^d)$.	Then 
	\[ F^g(y,U) = F(y,U) < \infty.\]
	So, there exists a sequence $u_k\in \mathcal{B}_{\e_k}(U)$ such that $u_k\to y$ in $L^p(\Omega;\R^d)$ and
	\[ \lim\limits_{k\to\infty} F_{\e_k}(u_k,U) = F(y,U). \]
	Let $\delta>0$, and then choose $U'$ such that $\supp \psi \subset U' \css U$ and $\abs{U\backslash U'}\leq \delta$.
	We now use a cut-off argument similarly as in the proof of Lemma \ref{lem:subadd}. Fix $N\in\N$ and define
	\[ U_j=\left\{x\in U \colon \dist(x,U') < \frac{j\dist(U', U^c)}{N} \right\}. \]
	Then choose the cut-off functions $\varphi_j \in C^{\infty}_c(U_{j+1}; [0,1])$ with $\varphi_j \equiv 1$ on $U_j$ and $\norm{\nabla \varphi_j}_{\infty} \le C N$ and set
	\[ \hat{g}_k(x) = g(a), \text{ if } a\in Q_{\e_k}(x)\cap \la_{\e_k} \text{ and} \]
	\[
	w_{n,j}(x)=
		\begin{dcases}
			\varphi_j(\bar{x}) u_k(x) + (1-\varphi_j(\bar{x})) \hat{g}_k(x), & \text{ if } x\in U^{\e_k},\\
			u_k(x) & \text{ otherwise.}
		\end{dcases}
	\]
	As in the proof of Lemma \ref{lem:subadd} we calculate
	\al{
		F_{\e_k}(w_{k,j},U) & \leq F_{\e_k}(u_k,U) + C(\norm{\nabla g}_\infty^p+1) \abs{U\backslash U'} + \underbrace{ \e_k^d \sum\limits_{\bar{x}\in (\la_{\e_k}'(W_j))^\circ} W_{\rm cell}(\bar{\nabla}w_{k,j}(\bar{x})) }_{:= S_{j,k}},
	}
	with $W_j=U_{j+2}\backslash\overline{U_{j-1}}$, estimate $S_{j,k}$ by averaging, choose $j(k)$ suitably and thus get
	\[
		\limsup\limits_{k\to\infty} F_{\e_k}(w_{k,j(k)},U) \leq F(y,U) + C \delta + \frac{C}{N-4}.
	\]
	Since we choose $j(k)\leq N-3$, we have $w_{k,j(k)}\in \mathcal{B}_{\e_k}(U,g)$ for any $k$ large enough. Furthermore, $w_{k,j(k)}\to y$ in $L^p(\Omega;\R^d)$ since $\psi$ has support in $U'$. Hence,
	\[
		\Gamma\text{-}\limsup\limits_{k\to\infty} F^g_{\e_k}(y,U) \leq F^g(y,U) + \delta C + \frac{C}{N-4}.
	\]
	Let $\delta\to 0$ and $N \to \infty$.

	In the general case $y|_U\in g + W_0^{1,p}(U;\R^d)$, take $y_l$ such that $y_l|_U \in g+ C_c^\infty(U;\R^d)$ and $y_l\to y$ in $W^{1,p}(U;\R^d)$ and in $L^p(\Omega;\R^d)$. We get
	\al{
		\Gamma\text{-}\limsup\limits_{k\to\infty} F^g_{\e_k}(y,U) &\leq \liminf\limits_{l\to\infty} (	\Gamma\text{-}\limsup\limits_{k\to\infty} F^g_{\e_k}(y_l,U) ) \\
		&\leq \liminf\limits_{l\to\infty} F^g(y_l,U) \\
		&= F^g(y,U)
	}
	by the lower semicontinuity of the $\Gamma$-$\limsup$ with respect to $L^p(\Omega;\R^d)$-convergence and continuity of $F^g(\cdot,U)$ with respect to $W^{1,p}(U;\R^d)$-convergence.
\ep

\subsection{The limiting minimum problem}

The following theorem is important in two ways. On the one hand we gain insight into the $\Gamma$-convergence result, on the other hand we will directly need it to get the homogenization result in Section \ref{sec:hom}.

\begin{thm} \label{thm:boundaryminconverge}
	Under the assumptions of Theorem \ref{thm:boundarygamma}, we have
	\[\min\limits_y F^g(y,U) = \lim\limits_{k\to\infty} (\inf\limits_y F_{\e_k}^g(y,U)).\]
	Furthermore, any sequence $y_k$ with equibounded energy	is pre-compact in $L^p(U;\R^d)$ and if we have a sequence satisfying
	\[\lim\limits_{k\to\infty} (\inf\limits_y F_{\e_k}^g(y,U)) = \lim\limits_{k\to\infty} F_{\e_k}^g(y_k,U),\]
	then every limit of a converging subsequence is a minimizer of $F^g(\cdot,U)$.
\end{thm}
\bp
	Fix $g,U$ and write $G_k(y)=F_{\e_k}^g(y,U)$, $G(y)=F^g(y,U)$. 
Let $y_k$ be a sequence with equibounded energy $G_k(y_k)$. By Assumption \ref{ass:coercive} and Proposition \ref{prop:gradients} we obtain that 
	\[ \int\limits_{U_{\e_{k}}} \abs{\nabla\tilde{y}_k}^p \,dx \leq C.\]
	Furthermore, using the boundary condition, we have
	\[ \int\limits_U \abs{\nabla\tilde{y}_k}^p \,dx \leq C.\]
	A Poincaré-type inequality involving the trace yields
	\al{
		\norm{\tilde{y}_k}_{W^{1,p}(U;\R^d)} &\leq C ( \norm{\nabla\tilde{y}_k}_{L^p(U;\R^d)} + \norm{T\tilde{y}_k}_{ L^p(\partial U;\R^d)} ) \\
		&\leq C + C \norm{g}_\infty \mathcal{H}^{d-1}(\partial U)^{\frac{1}{p}} \leq C 
	}
and so $\tilde{y}_{k_l} \to y$ in $L^p(U;\R^d)$ for some subsequence $k_l$ and some $y \in W^{1,p}(U; \R^d)$. Then, by Proposition \ref{prop:interpolbc}, $y_{k_l} \to y$ in $L^p(U;\R^d)$. 

Now from Theorem \ref{thm:boundarygamma} we infer that $G_k$ $\Gamma(L^p(\Omega;\R^d))$-converges to $G$. But then $G_k$ also $\Gamma(L^p(U;\R^d))$-converges to $G$. Here the existence of recovery sequences is immediate as $L^p(\Omega;\R^d)$- implies $L^p(U;\R^d)$-convergence. As for the $\liminf$-inequality, if $y_k \to y$ in $L^p(U; \R^d)$, where the energies $G_k(y_k)$ are without loss of generality assumed to be equibounded and, in particular, in $y_k\in\mathcal{B}_{\e_{k}}(U,g)$, we can extend the functions by defining $y_k(x) := g(\bar{x})$ and $y(x) := g(x)$ for $x\in \Omega\backslash U$ without changing their respective energies. Since then $y_k \to y$ in $L^p(\Omega; \R^d)$, we have indeed that $\liminf_{k \to \infty} G_k(y_k) \ge G(y)$. The remaining part of the proof now directly follows from Theorem \ref{thm:minconverge}. 
\ep

\section{Proof of the main results} \label{sec:hom}

\subsection{The $\Gamma$-convergence results}

To simplify notations, we define $P_h(x) = x + A (0,h)^d$ and $P_h = P_h(0)$.
First, we will prove the following lemma.
\begin{lem}\label{lem:homconvergence}
The limit
\[\lim\limits_{N\to\infty} \frac{1}{N^d} \inf\left\{ \sum\limits_{x\in (\la_1'(P_N))^\circ} W_{\rm cell}(\bar{\nabla}y(x)) \colon y \in \mathcal{B}_1(P_N,y_M) \right\}\]
exists for every $M\in\R^{d \times d}$.
\end{lem}
\bp
Let us define
\al{
	G(y,U) &= \sum\limits_{x\in (\la_1'(U))^\circ} W_{\rm cell}(\bar{\nabla}y(x))\quad \text{and} \\
	f_k(M) &= \frac{1}{k^d} \inf\left\{ G(y, P_k) \colon y \in \mathcal{B}_1(P_k,y_M) \right\}.
}
Fix $M\in\R^{d \times d}$ and let $k,n\in\N$ with $k > n$. Choose $v_n\in\mathcal{B}_1(P_n,y_M)$ such that
\[\frac{1}{n^d} G(v_n, P_n) \leq f_n(M) +\frac{1}{n}.\]
Now, we can define
\[
	u_k(x) =
	\begin{dcases}
		v_n(x - n \alpha) + n M\alpha &\text{ if } x\in P_n(n \alpha) \text{ for some } \alpha\in A\left\{ 0,1, \dots, \left[ \frac{k}{n}\right] -1 \right\}^d,\\
		M\bar{x} &\text{ otherwise.}
	\end{dcases}
\]
Since $v_n$ satisfies the boundary condition, $u_k$ is constant on every cell. Moreover, $u_k\in\mathcal{B}_1(P_k,y_M)$ and we can estimate
\al{
	f_k(M) &\leq \frac{1}{k^d} G (u_k,P_k) \\
	&\leq \frac{1}{k^d} \left( \left[ \frac{k}{n}\right]^d G (v_n,P_n) +  c (\abs{M}^p+1) \left( \#(\la_1'(P_k))^\circ - \left[ \frac{k}{n}\right]^d \#(\la_1'(P_n(\alpha n)))^\circ  \right) \right) \\
	&\leq \frac{1}{n^d} G (v_n,P_n) + \frac{c(\abs{M}^p+1)}{k^d} \left( \frac{\abs{P_k}- \abs{P_{n\left[ \frac{k}{n}\right]}}}{\abs{P_1}} + \left[ \frac{k}{n}\right]^d \left(n^d - (n-2)^d \right) \right) \\
	&\leq f_n(M) +\frac{1}{n} + \frac{c(\abs{M}^p+1)}{k^d} \left( k^d - \left( n\left[ \frac{k}{n}\right]\right)^d + k^d \left(1 - \left(1-\frac{2}{n}\right)^d \right) \right)  \\
	&\leq f_n(M) +\frac{1}{n} + c(\abs{M}^p+1) \left( 1 - \left( 1-\frac{n}{k} \right)^d + 1 - \left(1-\frac{2}{n}\right)^d \right).
}
Thus, for every $n\in\N$,
\[ \limsup\limits_{k\to\infty} f_k(M) \leq f_n(M) +\frac{1}{n} + c(\abs{M}^p+1) \left( 1 - \left(1-\frac{2}{n}\right)^d \right), \]
hence,
\[\limsup\limits_{k\to\infty} f_k(M) \leq \liminf\limits_{n\to\infty} f_n(M).\]
\ep

Now, we can prove our first main theorem.
\bp[Proof of Theorem \ref{thm:homogenization}]
We will first show that $F_{\e}(\cdot, \Omega)$ $\Gamma(L^p(\Omega; \R^d))$-converges to $F$.
According to Lemma \ref{lem:homconvergence}, $W_{\rm cont}$ is well-defined. 
By the Urysohn property of $\Gamma$-convergence in Proposition \ref{prop:gammaproperties}, it is enough to show that, for any sequence $\e_n\to 0$, the function $f$ of Theorem \ref{thm:compactness} equals $W_{\rm cont}$. Fix such a sequence, the subsequence $\e_k$ and the associated $f$. Since $f$ is quasiconvex, we have for every $M\in\R^{d \times d}$ and $U\in\mathcal{A}_L(\Omega)$
\al{
	f(M) &= \frac{1}{\abs{U}} \min\left\{ \int\limits_U f(\nabla y (x)) \,dx \colon y-y_M \in W^{1,p}_0(U;\R^d) \right\} \\
	&= \frac{1}{\abs{U}} \min\left\{ F(y,U) \colon y-y_M \in W^{1,p}_0(U;\R^d) \right\} .
}
If we restrict $y_M$ to a ball that contains some neighborhood of $\Omega$, we can extend it to a function in $W^{1,\infty}(\R^d;\R^d)\cap C(\R^d;\R^d)$, so $y_M$ is admissible as a boundary condition in Theorem \ref{thm:boundarygamma} and we get the $\Gamma$-convergence result with boundary condition. Hence by Theorem \ref{thm:boundaryminconverge}, for $h_0>0$ and $x_0\in \R^d$ such that $P_{h_0}(x_0)\css \Omega$
\al{
	f(M) &= \frac{1}{\abs{P_{h_0}(x_0)}} \lim\limits_{k\to\infty} \left(\inf\left\{ F_{\e_k}(y,P_{h_0}(x_0)) \colon y \in \mathcal{B}_{\e_k}(P_{h_0}(x_0),y_M) \right\}\right) \\
	&=\frac{1}{\abs{P_{h_0}(x_0)}} \lim\limits_{k\to\infty} \inf\left\{ \e_k^d \sum\limits_{x\in (\la_{\e_k}'(P_{h_0}(x_0)))^\circ} W_{\rm cell}(\bar{\nabla}y(x)) \colon y \in \mathcal{B}_{\e_k}(P_{h_0}(x_0),y_M) \right\}
}
It is easy to see, that we can always find $h_k > 0$ and $x_k\in \la_{\e_k}$ such that
\[P_{h_k}(x_k) = \left(\bigcup\limits_{x\in P_{h_0}(x_0)} Q_{\e_k}(x)\right)^\circ .\]
We then know $P_{h_k}(x_k) \subset \Omega$ for all $k$ large enough, $\abs{x_0-x_k}\leq \diam{Q_{\e_k}} = \e_k \diam{Q_1}$, $h_0\leq h_k \leq h_0 + 2 \e_k$ and, that there are $N_k\in \N$ satisfying $h_k = N_k \e_k$. Furthermore,
\al{
	\la_{\e_k}'(P_{h_0}(x_0)) &= \la_{\e_k}'(P_{h_k}(x_k)) \quad \text{and}\\
	(\la_{\e_k}'(P_{h_0}(x_0)))^\circ &= (\la_{\e_k}'(P_{h_k}(x_k)))^\circ.
}
Hence, $\mathcal{B}_{\e_k}(P_{h_0}(x_0),y_M)$ and $\mathcal{B}_{\e_k}(P_{h_k}(x_k),y_M)$ are equal up to extending the functions in $\mathcal{B}_{\e_k}(P_{h_0}(x_0),y_M)$ constant on cells that intersect $P_{h_k}(x_k) \setminus P_{h_0}(x_0)$. It follows that
\al{
	f(M) &=\frac{1}{\abs{P_1}} \lim\limits_{k\to\infty} \frac{1}{N_k^d} \frac{h_k^d}{h_0^d} \inf\left\{ \sum\limits_{x\in (\la_{\e_k}'(P_{h_k}(x_k)))^\circ} W_{\rm cell}(\bar{\nabla}y(x)) \colon y \in \mathcal{B}_{\e_k}(P_{h_k}(x_k),y_M) \right\}\\
	&=\frac{1}{\abs{P_1}} \lim\limits_{k\to\infty} \frac{1}{N_k^d} \inf\left\{ \sum\limits_{x\in (\la_{\e_k}'(P_{h_k}))^\circ} W_{\rm cell}(\bar{\nabla}y(x)) \colon y \in \mathcal{B}_{\e_k}(P_{h_k},y_M) \right\},
}
where we used, that $y \in \mathcal{B}_{\e_k}(P_{h_k}(x_k),y_M)$, if and only if $y(\cdot + x_k) - Mx_k \in \mathcal{B}_{\e_k}(P_{h_k},y_M)$ and that the discrete gradient of $y$ at a point $x$ equals the discrete gradient of $y(\cdot + x_k) - Mx_k$ at $x-x_k$.
In a similar way $y \in \mathcal{B}_{\e_k}(P_{h_k},y_M)$ if and only if $y' \in \mathcal{B}_1(P_{N_k},y_M)$ and $\bar{\nabla}y'(x) = \bar{\nabla}y(\e_k x)$, where $y'(x) = \frac{1}{\e_k} y(\e_k x)$. Hence,
\al{
	f(M) &=\frac{1}{\abs{P_1}} \lim\limits_{k\to\infty} \frac{1}{N_k^d} \inf\left\{ \sum\limits_{x\in (\la_1'(P_{N_k}))^\circ} W_{\rm cell}(\bar{\nabla}y(x)) \colon y \in \mathcal{B}_1(P_{N_k},y_M) \right\}\\
	&= W_{\rm cont}(M).
}

In order to prove that also $F_{\e}(\cdot, \Omega)$ $\Gamma(L^p_{loc}(\Omega; \R^d)/\R)$-converges to $F$, we only need to verify the $\liminf$-inequality as the existence of recovery sequences immediately follows from the first part of the proof since convergence in $L^p(\Omega; \R^d)$ implies convergence in $L^p_{loc}(\Omega; \R^d)/\R$. But if $\e_n \to 0$ and $y_{\e_n} \to y$ in $L^p_{loc}(\Omega; \R^d)/\R$, then there exist $c_n \in \R$ such that, for every $U \in \mathcal{A}_L(\Omega)$ with $U \css \Omega$, $y_{\e_n} - c_n \to y$ in $L^p(U; \R^d)$, so that by the previous result 
\[ \liminf\limits_{n \to \infty} F_{\e_n}(y_{\e_n}, \Omega) 
= \liminf\limits_{n \to \infty} F_{\e_n}(y_{\e_n} -c_n, \Omega) 
\ge \liminf\limits_{n \to \infty} F_{\e_n}(y_{\e_n}- c_n, U) 
\ge F(y, U). \] 
Without loss of generality we may assume that $\liminf\limits_{k \to \infty} F_{\e_k}(y_{\e_k}, \Omega) < \infty$. Since for any $V \in \mathcal{A}(\Omega)$ with $V \css \Omega$ there exists $U \in \mathcal{A}_L(\Omega)$ with $V \subset U \css \Omega$, we then deduce from Lemma \ref{lem:gammaliminf} that $y \in W^{1,p}(V; \R^d)$ with $\norm{y}_{W^{1,p}(V; \R^d)}$ bounded uniformly in $V \in \mathcal{A}$ with $V \css \Omega$, hence $y \in W^{1,p}(\Omega; \R^d)$. Then invoking Lemma \ref{lem:limitregular} and passing to the supremum over $U \in \mathcal{A}_L(\Omega)$ in the above inequality yields 
\[  \liminf\limits_{k \to \infty} F_{\e_k}(y_{\e_k}, \Omega) \ge F(y, \Omega). \]
\ep

\bp[Proof of Theorem \ref{thm:homogenization-bdry}]
Theorem \ref{thm:homogenization-bdry} is a direct consequence of Theorem \ref{thm:boundarygamma} and Theorem \ref{thm:homogenization}, where the limiting energy density $f$ has been identified as $W_{\rm cont}$.
\ep

\bp[Proof of Theorem \ref{thm:bd-En-seq}]
Suppose $y_k$ is a sequence with equibounded energies $F_{\e_k}(y_k)$. By Proposition \ref{prop:gradients} and the growth assumptions on $W_{\rm cell}$, for every $U \in \mathcal{A}(\Omega)$ with $U \css \Omega$ we have 
\[ \int_U \abs{\nabla \tilde{y}_k}^p 
\le C F_{\e_k}(y_k) + C|\Omega| \]
uniformly bounded for sufficiently large $k$. Choose $U_0 \in \mathcal{A}_L(\Omega)$ connected and with $\emptyset \neq U_0 \css \Omega$. As $U_0$ is connected, by Poincaré's inequality we find $c_k \in \R$ such that $\tilde{y}_k - c_k$ is pre-compact in $L^p(U_0; \R^d)$. But then indeed for any connected $U \in \mathcal{A}_L$ with $U_0 \subset U \css \Omega$ the Poincar{\'e} inequality 
\[ \norm{\tilde{y}_k - c_k}_{W^{1,p}(U; \R^d)} 
\le C \norm{\nabla \tilde{y}_k}_{L^p(U; \R^d)} + \norm{\tilde{y}_k - c_k}_{L^p(U_0; \R^d)} \] 
yields that $\tilde{y}_k - c_k$ is pre-compact in $L^p(U; \R^d)$. Exhausting $\Omega$ with a countable number of such domains and passing to a diagonal sequence, we find a subsequence $y_{k_n}$ such that $\tilde{y}_{k_n} - c_{k_n}$ converges in $L^p_{loc}(\Omega; \R^d)$. By Remark \ref{rem:Lploc-conv} we finally obtain that $y_{k_n} - c_{k_n}$ converges in $L^p_{loc}(\Omega; \R^d)$. 
\ep

\bp[Proof of Theorem \ref{thm:bd-En-seq-bdry}]
This is immediate from Theorem \ref{thm:boundaryminconverge}. \ep

\bp[Proof of Corollary \ref{cor:conv-almost-min}]
This is a direct consequence of Theorems \ref{thm:homogenization}, \ref{thm:bd-En-seq}, \ref{thm:homogenization-bdry}, \ref{thm:bd-En-seq-bdry} and \ref{thm:minconverge}. 
\ep

\bp[Proof of Theorem \ref{thm:CB}]

If in addition to Assumptions \ref{ass:coercive} and \ref{ass:growth} Assumption \ref{ass:CB} holds true, we can apply \cite[Theorem 4.2]{CDKM} with $\Lambda = (\la_1'(P_{N_k}))^\circ$. It is easy to see that the boundary of $\Lambda$ as defined in \cite{CDKM} equals $\partial\la_1(P_{N_k})\cup \la_1 \backslash \la_1'(P_{N_k})$, but of course the second part does not change anything. This shows that there is a neighborhood ${\cal U}$ of $SO(d)$, such that for every $M\in {\cal U}$
\al{
	W_{\rm cont}(M) 
	&=\frac{1}{\abs{P_1}} \lim\limits_{k\to\infty} \frac{1}{N_k^d} \inf\left\{ \sum\limits_{x\in (\la_1'(P_{N_k}))^\circ} W_{\rm cell}(\bar{\nabla}y(x)) \colon y \in \mathcal{B}_1(P_{N_k},y_M) \right\}\\
	&= \frac{1}{|\det A|} \lim\limits_{k\to\infty} \frac{1}{N_k^d} \sum\limits_{x\in (\la_1'(P_{N_k}))^\circ} W_{\rm cell}(MZ) \\
	&= \frac{1}{|\det A|} W_{\rm cell}(MZ) \\
	&= W_{\rm CB}(M).
}
\ep

\subsection{Approximation of general continuum densities}

Next we prove Propositions \ref{prop:all-qc} and \ref{prop:all-quad}. 
\bp[Proof of Proposition \ref{prop:all-qc}] 
At variance with our previous decomposition procedure, we now choose any simplicial decomposition $\mathcal{S}$ of the cell $A[-\frac{1}{2}, \frac{1}{2})^d$ into $d$-simplices all of whose corners lie in $A\{-\frac{1}{2}, \frac{1}{2}\}^d$. For $F = (f_1, \ldots, f_{2^d}) \in \R^{d \times 2^d}$ we then interpolate the mapping 
\[ A\left\{-\frac{1}{2}, \frac{1}{2}\right\}^d \to \R^d, \qquad x_i \mapsto f_i \]
affine on each simplex in order to obtain 
\[ u_F : A\left[-\frac{1}{2}, \frac{1}{2}\right)^d \to \R^d. \]
Then $W_{\rm cell}$ is defined by 
\[ W_{\rm cell}(F) := \int_{A[-\frac{1}{2}, \frac{1}{2})^d} V(\nabla u_F) \, dx. \]
As every corner $z_{i_0},\ldots, z_{i_d}$ of $S\in\mathcal{S}$ lies in $A\{-\frac{1}{2}, \frac{1}{2}\}^d$, we have \[  c \sum\limits_{j=1}^d \abs{f_{i_j}-f_{i_0}} 
\leq \abs{\nabla u_F} 
\leq C\abs{F} \]
on $S$. Thus, $\norm{F}=\max\limits_{x \in A[-\frac{1}{2}, \frac{1}{2})^d} \abs{\nabla u_F(x)}$ is a norm on $V_0$ and we calculate
\al{ W_{\rm cell}(F) &\geq \int_{A[-\frac{1}{2}, \frac{1}{2})^d} c\abs{\nabla u_F}^p - c' \, dx \\
&\geq c \norm{F}^p -c' \\
&\geq c \abs{F}^p -c',
}
and on the other hand
\[  W_{\rm cell}(F) \leq C(\abs{F}^p + 1).\]
This means $W_{\rm cell}$ satisfies Assumptions \ref{ass:coercive} and \ref{ass:growth}. From Theorem \ref{thm:homogenization} we then deduce that 
\al{
	W_{\rm cont}(M) 
	&=\frac{1}{\abs{\det A}} \lim\limits_{N\to\infty} \frac{1}{N^d} \inf\left\{ \sum\limits_{x\in (\la_1'(P_{N}))^\circ} W_{\rm cell}(\bar{\nabla}y(x)) \colon y \in \mathcal{B}_1(P_N,y_M) \right\}\\
	&=\frac{1}{\abs{\det A}} \lim\limits_{N\to\infty} \frac{1}{N^d} \inf\left\{ \int_{(P_N)_1} V( \nabla u_{\bar{\nabla}y(x)})  \colon y \in \mathcal{B}_1(P_N,y_M) \right\} \\
	&=\frac{1}{\abs{\det A}} \lim\limits_{N\to\infty} \frac{1}{N^d} \abs{\det(A)} (N-2)^d V(M) \\
	&= V(M) 
} 
due to the quasiconvexity of $V$. 
\ep

\bp[Proof of Proposition \ref{prop:all-quad}] 
Any $F \in V_0$ can be decomposed orthogonally as $F = F'Z + F''$ with unique $F' \in \R^{d \times d}$ and $F'' \in (\R^{d \times d}Z)^{\perp}$. Set
\[ W_{\rm cell}(F) 
= \abs{\det A} Q\left(\sqrt{(F')^T F'} - \Id\right) + |F''|^2 + \chi(F), \] 
where $\chi$ is any frame indifferent function satisfying Assumptions \ref{ass:coercive} and \ref{ass:growth} with $p \ge d$ which is non-negative, vanishes near $\bar{SO}(d)$ and is bounded from below by a positive constant on $\bar{O}(d) \setminus \bar{SO}(d)$, $\bar{O}(d) = O(d)Z$. Then also $W_{\rm cell}$ satisfies Assumptions \ref{ass:coercive} and \ref{ass:growth} with the same $p$. Noting that, for $M \in \R^{d \times d}$, $(MF)' = M F'$ and $(MF)'' = MF''$, it is not hard to verify that $W_{\rm cell}$ also satisfies Assumption \ref{ass:CB} with 
\[ D^2 W_{\rm cell}(Z)(F, F) 
= 2 \abs{\det A} Q\left(\frac{(F')^T + F'}{2}\right) + 2|F''|^2.\]   
But then 
\[ \frac{1}{2} D^2 W_{\rm CB}(\Id)(M,M) 
= \frac{1}{2 \abs{\det A}} D^2 W_{\rm cell}(Z)(MZ, MZ) 
= Q\left(\frac{M^T + M}{2}\right) 
= Q(M).\]   
\ep

\subsection{Extension to finite-range energies}

We briefly comment on more general long-range interactions. Suppose $\Lambda = \{z_1, \ldots, z_{2^d}, \ldots, z_N\} \subset \la'$ is any fixed finite set, where $z_1, \ldots, z_{2^d}$ still denote $A\{-\frac{1}{2}, \frac{1}{2}\}^d$. For $y \in \mathcal{B}_{\e}(\Omega)$ we define $y_i = y(\bar{x} + \e z_i)$. With $\bar{x}$ and $\bar{y}$ as before, i.e., only depending on $y_1, \ldots y_{2^d}$, let now 
\[\bar{\nabla}y(x) = \frac{1}{\e} (y_1-\bar{y}, \dots, y_{N}-\bar{y}) \in \R^{d \times N}. \]

The lattice interior $(\la_\e'(U))^\circ$ and boundary $\partial\la_\e'(U)$ now have to be shrunk respectively enlarged to a whole boundary layer, according to the maximal interaction length in $\Lambda$. Assumptions \ref{ass:coercive} and \ref{ass:growth} are then replaced by the estimate 
\[ c\abs{F'}^p-c' \leq W_{\rm super-cell}(F) \leq c''(\abs{F}^p+1)\]
for constants $c, c', c'' > 0$ and all $F\in \R^{d \times N}$ which satisfy $F'\in V_0$, where $F' \in \R^{d\times 2^d}$ denotes the left $d \times 2^d$ submatrix of $F$. Note that the lower bound in particular allows for arbitrarily weak long range interactions. As the interpolation we used only depends on the $d \times 2^d$ values of the corresponding lattice cell, this implies that we get the standard estimates for the gradients in Proposition \ref{prop:gradients} only on this part of the discrete gradient.

It is important that the interaction range is bounded by $C \e$, so that, e.g., Lemma \ref{lem:growth} and its proof still work.
In the estimates of the error $S_{j,n}$ in, e.g., Lemma \ref{lem:subadd}, it is important that $\e^{-1}\abs{u(\bar{x}+\e z_i)-u(\bar{x})}$ and thus the discrete gradient can be bounded by a fixed finite sum of smaller $d \times 2^d$ discrete gradients of some cells near $x$. Hence, we still have the estimate
\[  \e_n^d \sum\limits_{x\in (\la_{\e_n}'(U))^\circ} \abs{\bar{\nabla}u(x)}^p
		\leq C \int\limits_{U}  \abs{\nabla\tilde{u}(x)}^p \,dx \]  
Note that according to our enlarging of the lattice boundaries, also the cell formula for the limit density will now involve a sequence of minimizing problems with affine boundary conditions on a boundary layer.

 We finally remark that the statement on the applicability of the Cauchy-Born rule translates naturally, as the main ingredient does, see \cite[Theorem 5.1]{CDKM}.

\subsection{Extension to multi-lattices}
It is also possible to generalize these results to certain non-Bravais lattices, namely to multi-lattices of the form $\la\cup (s_1+\la)\cup \dots \cup (s_m+\la)$, in the following way:
We still consider $\la$ to be our main lattice. But now we have $m$ additional atoms in each cell, which we describe by the `internal variable' $s(x)\in\R^{d\times m}$, such that $\e s_{\cdot j}$ describes the distance of the $j$-th atom to the midpoint of the cell. Of course $s$ can be identified with a function, that is constant on every interior cell and is $0$ outside and thus lies in some $L^q(\Omega;\R^{d\times m})$, $1<q<\infty$. The new cell energy depends on $md$ additional variables and we now consider the growth condition
\[c(\abs{M'}^p + \abs{s}^q) - c' \leq W_{\rm super-cell}(M,s)\leq c''(\abs{M}^p + \abs{s}^q + 1)\]
for $M\in \R^{d\times N}$ and $s\in \R^{d\times m}$. It is now natural to have a $\Gamma$-convergence result with respect to strong-$L^p$-convergence in the first and weak-$L^q$-convergence in the second component.
As we will see in a moment, it turns out that we have to consider a combined boundary value and mean value problem. For this we define $\mathcal{B}_\e(U,g,s_0)$ to consist of all pairs $(y,s)$, such that $y\in\mathcal{B}_\e(U,g)$ and $s\in L^q(\Omega;\R^{d\times m})$ is constant on every interior cell of $U$, is $0$ outside and has mean value $s_0$ on the union of interior cells of $U$.

In analogy to Theorem \ref{thm:homogenization}, we now have Theorem \ref{thm:homogenization-multicell}. The proof of this theorem is similar to the proof of Theorem \ref{thm:homogenization}. But there are several things that need to be addressed:

First of all, the weak topology on $L^q$ is not given by a metric. But, as discussed in \cite{dalmasogamma} in detail, this is not a big problem, since our functionals are equicoercive and the dual of $L^q$ is separable. In particular, we can describe $\Gamma$-convergence by sequences and the compactness and the Urysohn property are still true. 
Next, we need an advanced version of our integral representation result:
\begin{thm}
	Let $1\leq p,q < \infty$ and let $F\colon W^{1,p}(\Omega;\R^d)\times L^q(\Omega;\R^{d\times m})\times\mathcal{A}(\Omega) \to [0,\infty]$ satisfy the following conditions:
	\begin{enumerate}[(i)]
		\item (locality) $F(y,s,U) = F(v,t,U)$, if $y(x) = v(x)$ and $s(x)=t(x)$ for a.e.\ $x\in U$;
		\item (measure property) $F(y,s,\cdot)$ is the restriction of a Borel measure to $\mathcal{A}(\Omega)$;
		\item (growth condition) there exists $c>0$ such that
			\[ F(y,s,U)\leq c \int\limits_U \abs{\nabla y(x)}^p +\abs{s}^q + 1 \,dx;\]
		\item (translation invariance in y) $F(y,s,U) = F(y+a,s,U)$ for every $a\in\R^d$ ;
		\item (lower semicontinuity) $F(\cdot,\cdot,U)$ is sequentially lower semicontinuous with respect to weak convergence in $W^{1,p}(\Omega;\R^d)$ in the first and weak convergence in $L^q(\Omega;\R^{d\times m})$ in the second component;
		\item (translation invariance in x) With $y_M(x)=M x$ and $s(x)=s_0$ we have
		\[ F(y_M,s,B_r(x))=F(y_M,s,B_r(x'))\]
		for every $M\in\R^{d\times d}$, $s_0\in \R^{d \times m}$, $x,x'\in \Omega$ and $r>0$ such that $B_r(x),B_r(x')\subset\Omega$.
	\end{enumerate}
	Then there exists a continuous $f \colon \R^{d\times d}\times \R^{d \times m} \to [0,\infty)$ such that
	\[
		0\leq f(M,s) \leq C(1+\abs{M}^p+\abs{s}^q)\]
	for every $M\in\R^{d\times d}$, $s\in\R^{d\times m}$ and
	\[
		F(y,s,U) = \int\limits_U f(\nabla y(x),s(x))\,dx
	\]
	for every $y\in W^{1,p}(\Omega;\R^d)$, $s\in L^q(\Omega;\R^{d\times m})$ and $U\in\mathcal{A}(\Omega).$
\end{thm}
The proof in \cite{braidesdefranceschi} for the pure Sobolev version of this theorem readily applies to this more general statement. (Note that continuity of $f$ then follows from seperate convexity.) 
Most of the lemmata then translate naturally. We just want to comment on some details in Lemma \ref{lem:subadd}. The recovery sequences now contain additionally some $t_n \wto s$, $r_n \wto s$ corresponding to $U$,$V$ respectively. We define \[q_{n,j} = \chi_{U_j} (\bar{x}) t_n(x) + (1 - \chi_{U_j} (\bar{x})) r_n(x),\] and then choose $j(n)$ as before to define $s_n=q_{n,j(n)}$. The only part that is not immediately clear now, is the convergence $s_n \wto s$. To prove this, let $\varphi\in L^{q'}(\Omega ; \R^{d\times m})$. We now split $\varphi$ into several parts we can control
\[\varphi = \psi_n +  \varphi\chi_{U'} + \varphi\chi_{\Omega \backslash \overline{U_N}} + \sum\limits_{j=0}^{N-1} \varphi \chi_{(U_{j+1}\backslash \overline{U_j})_{\e_n}}.\]
Here the $\psi_n$ contain all the remaining parts. We see that $\psi_n \to 0$ strongly in $L^{q'}$ as long as $\abs{\partial U_j}=0$ for every $j$, so this is true up to changing the sets $U_j$ a little bit. But then we also have
\[\varphi\chi_{(U_{j+1}\backslash \overline{U_j})_{\e_n}} \to \varphi\chi_{U_{j+1}\backslash U_j}\]
strongly in $L^{q'}$. The advantage is now that on each set $(U_{j+1}\backslash \overline{U_j})_{\e_n}$ we have either $s_n = t_n$ or $s_n = r_n$, possibly changing with $n$. But in both cases we have weak convergence to $s$, hence
\[\int\limits_\Omega s_n(x) \varphi(x) \chi_{(U_{j+1}\backslash \overline{U_j})_{\e_n}}(x)\,dx \to \int\limits_\Omega s(x) \varphi(x) \chi_{U_{j+1}\backslash U_j}(x)\,dx.\] 
And, putting it all together, we get
\[\int\limits_\Omega s_n(x) \varphi(x)\,dx \to \int\limits_\Omega s(x) \varphi(x)\,dx. \]

Another important step is to adjust Theorem \ref{thm:boundarygamma} and Theorem \ref{thm:boundaryminconverge}, so that additionally to the boundary values for $y$ we have a fixed mean value for $s$, i.e., we consider $\mathcal{B}_\e(U,g,s_0)$ instead of $\mathcal{B}_\e(U,g)$ in the discrete setting and add the constraint \[\fint\limits_U s(x)\,dx = s_0\]
in the continuum setting. To get the $\liminf$-inequality just notice that for $s_k\wto s$ with $(y_k,s_k)\in\mathcal{B}_\e(U,g,s_0)$ we have
\al{\fint\limits_U s(x)\,dx&= \lim\limits_{k\to\infty}\fint\limits_U s_k(x)\,dx\\
&=\lim\limits_{k\to\infty}\frac{\abs{U_{\e_k}}}{\abs{U}}\fint\limits_{U_{\e_k}} s_k(x)\,dx\\
&=\lim\limits_{k\to\infty}\frac{\abs{U_{\e_k}}}{\abs{U}} s_0 = s_0.
}
The $\limsup$-inequality is a little more subtle. We have a function $s\in L^q(\Omega;\R^{d\times m})$ with \[\fint\limits_U s(x)\,dx = s_0\]
and a recovery sequence without this mean value $s_k\wto s$. Let us write
\[\fint\limits_{U_{\e_k}} s_k(x)\,dx + \xi_k = s_0, \]
so that $\xi_k \to 0$. We now adjust the $s_k$ adequately.
Define
\[t_k(x) = s_k(x) + \xi_k\frac{\abs{U_{\e_k}}}{\abs{V_k}}\chi_{V_k}.\]
If $V_k$ is a union of cells with some distance to the boundary of $U$, then, for $k$ large enough, the $t_k$ are admissible functions and do not interact with the adjustments on $y$. We have to make sure, that $t_k\wto s$ and \[\limsup\limits_{k\to\infty} F_{\e_k}(u_k,t_k,U) \leq \limsup\limits_{k\to\infty} F_{\e_k}(u_k,s_k,U).\]
The weak convergence is true, if $\abs{V_k}\to 0$ and $\abs{V_k}\geq c \xi_k$ for some $c>0$. For the second estimate, we have to choose the $V_k$ a little more carefully to avoid concentration of the energy. Choose sequences $\eta_k \to 0$, $L_k \to \infty$ such that $\eta_k \ge c\xi_k$, $\frac{\eta_k}{\e_k^d} \to \infty$ and $L_k \eta_k \to 0$. Then take $L_k\in\N$ disjoints sets $W_{k,l}\subset U$, that are unions of cells, such that $\abs{W_{k,l}}$ is independent of $l$ and is roughly equal to $\eta_k$, which means
\[c\eta_k \leq \abs{W_{k,l}} \leq C \eta_k,\]
with $C,c>0$ independent of $k$ and $l$. 
This is possible as $\frac{\eta_k}{\e_k^d} \to \infty$ and $L_k \eta_k \to 0$. Then, we can choose $l(k)$ and set $V_k = W_{k,l(k)}$, such that
\al{
		&\int\limits_{V_k} W_{\rm cell}(\bar{\nabla}u_k(x),s_k(x)) + W_{\rm cell}\big(\bar{\nabla}u_k(x),s_k(x)+\xi_k\frac{\abs{U_{\e_k}}}{\abs{V_k}}\big)\, dx \\
		&\leq \frac{1}{L_k} \sum\limits_{l=1}^{L_k} \int\limits_{W_{k,l}} W_{\rm cell}\big(\bar{\nabla}u_k(x),s_k(x)) + W_{\rm cell}(\bar{\nabla}u_k(x),s_k(x)+\xi_k\frac{\abs{U_{\e_k}}}{\abs{W_{k,l}}}\big)\, dx \\
		&\leq \frac{1}{L_k} C,
	}
due to the growth condition. So the error goes to zero with $L_k\to\infty$. 
The rest of the proof translates naturally. The most important observation is the equality
\al{
	f(M,s_0) 
	&= \frac{1}{\abs{U}} \min\Bigg\{ \int\limits_U f(\nabla y (x),s(x)) \,dx \colon y-y_M \in W^{1,p}_0(U;\R^d),\\
	&\qquad\qquad\qquad\qquad\qquad s\in L^q(U;\R^{d\times m}), \fint\limits_U s(x)\,dx = s_0 \Bigg\} \\
	&= \frac{1}{\abs{U}} \min\Bigg\{  F(y,s,U) \colon y-y_M \in W^{1,p}_0(U;\R^d),\\
	&\qquad\qquad\qquad\qquad\qquad s\in L^q(U;\R^{d\times m}), \fint\limits_U s(x)\,dx = s_0 \Bigg\},
}
which is of course a consequence of the lower semicontinuity properties.

\bp[Proof of Theorem \ref{thm:mins}]
	Fix $y\in W^{1,p}(\Omega;\R^d)$ and without loss of generality fix a version of $y$ that is finite everywhere.

	Due to the growth condition and the continuity, we know that the infimum in
	\[\inf\limits_{s\in \R^{d\times m}} W_{\rm cont}(M,s)\]
	is actually a minimum for arbitrary $M$ and that the function
	\[M\mapsto \min\limits_{s\in \R^{d\times m}} W_{\rm cont}(M,s)\]
	is continuous.
	Obviously, we always have the inequality	
	\[ \int\limits_\Omega W_{\rm cont}(\nabla y(x),s(x))\, dx \geq \int\limits_\Omega \min\limits_{s\in \R^{d\times m}} W_{\rm cont}(\nabla y (x),s)\, dx.\]
	We now want to show, that there always exists an $L^q$-function $s$ where this is an equality. The idea is of course to choose $s(x)$ as a minimizer of  $s \mapsto W_{\rm cont}(\nabla y (x),s)$. The key point is to ensure measurability. We will do this by using the theory of measurable multifunctions as developed, e.g., in \cite{fonsecaleonilp}. Define 
\[  \Theta(M)=\{s\in\R^{d\times m} \colon W_{\rm cont}(M,s) = \min\limits_{t\in \R^{d\times m}} W_{\rm cont}(M,t) \}
\]
and set $\Gamma (x) = \Theta ( \nabla y(x))$. Due to the continuity and the growth of $W_{\rm cont}$, the set $\Gamma(x)$ is closed and non-empty for every $x\in\Omega$, hence $\Gamma \colon \Omega \to \mathcal{P}(\R^{d\times m})$ is a closed-valued multifunction.

Next, we want to show that $\Gamma$ is measurable, in the sense that
\[ \Gamma^-(C)=\{x\in\Omega \colon \Gamma(x)\cap C \neq \emptyset\} \]
is Lebesgue-measurable for every closed set $C\subset \R^{d\times m}$. To this end, we will first show that $\Theta^- (C)$ is closed. Let $M_n\in\Theta^- (C)$, $M_n\to M$ and choose $s_n \in \Theta(M_n)\cap C$. Using the growth of $W_{\rm cont}$ and since the $M_n$ are bounded, the $s_n$ are also bounded. So for some subsequence we have $s_{n_k}\to s$ and $s\in C$. Furthermore,
\al{ W_{\rm cont}(M,s) &= \lim\limits_{k\to \infty} W_{\rm cont}(M_{n_k},s_{n_k})\\
&= \lim\limits_{k\to \infty} \min\limits_{t\in\R^{d\times m}} W_{\rm cont}(M_{n_k},t)\\
&= \min\limits_{t\in\R^{d\times m}} W_{\rm cont}(M,t).
}
This proves $s\in \Theta(M)\cap C$, so $\Theta^-(C)$ is closed and $\Gamma^-(C)=(\nabla y)^{-1} ( \Theta^-(C))$ is Lebesgue-measurable. Now, we can apply \cite[Thm.~6.5]{fonsecaleonilp}, to get a measurable $s\colon \Omega \to \R^{d\times m}$, with \[W_{\rm cont}(\nabla y (x),s(x)) = \min\limits_{t\in \R^{d\times m}} W_{\rm cont}(\nabla y (x),t)\] 
and $s\in L^q(\Omega; \R^{d\times m})$, since
\al{
	\int\limits_\Omega \abs{s(x)}^q\,dx &\leq C \int\limits_\Omega W_{\rm cont}(\nabla y(x),s(x)) + 1\,dx \\
	&= C \int\limits_\Omega \min\limits_{s\in\R^{d \times m}} W_{\rm cont}(\nabla y(x),s) + 1\,dx \\
	&\leq C \int\limits_\Omega \min\limits_{s\in\R^{d \times m}} \abs{\nabla y(x)}^p +\abs{s}^q + 1\,dx \\
	&\leq C \int\limits_\Omega \abs{\nabla y(x)}^p + 1\,dx.
}

It remains to justify the $\Gamma$-convergence result for $F^{s-\min}_\e$. Suppose $y_k \to y \in W^{1, p}(\Omega; \R^d)$ strongly in $L^p(\Omega; \R^d)$. Choose $s_k \in L^q(\Omega; \R^{d \times m})$ with $F^{s-\min}_{\e_k}(y_k, \Omega) \le F_{\e_k}(y_k, s_k, \Omega) + k^{-1}$. Without loss of generality assuming that $F^{s-\min}_\e(y_k, \Omega)$ is bounded, by passing to a subsequence (not relabeled) we may assume that $s_k \wto s_0$ in $L^q$. But then Theorem \ref{thm:homogenization-multicell} shows that 
$$ \liminf_{k \to \infty} F^{s-\min}_{\e_k}(y_k, \Omega) 
   = \liminf_{k \to \infty} F_{\e_k}(y_k, s_k, \Omega) 
   \ge F(y, s_0, \Omega) 
   \ge F^{s-\min}(y, \Omega) $$ 
by the first part of the proof. On the other hand, if $y \in W^{1, p}(\Omega; \R^d)$ is given, choose $s\in L^q(\Omega; \R^{d \times m})$ according to the first part of the proof such that $F^{s-\min}(y, \Omega) = F(y, s, \Omega)$. Then if $(y_k, s_k)$ is a recovery sequence for $(y, s)$ from Theorem \ref{thm:homogenization-multicell}, we obtain 
\[ \limsup_{k \to \infty} F^{s-\min}_{\e_k}(y_k, \Omega) 
   \le \limsup_{k \to \infty} F_{\e_k}(y_k, s_k, \Omega)
   = F(y, s, \Omega) 
   = F^{s-\min}(y, \Omega).\]
\ep

\bibliographystyle{alpha}
\bibliography{disc-cont-nonlin}

\end{document}